\documentclass[11pt]{article} 
\usepackage[a4paper, margin=2.0cm]{geometry}

\usepackage{amssymb,amsmath,amsthm}
\usepackage{bm}

\usepackage[authoryear,round]{natbib}
\setcitestyle{authoryear,round}

\usepackage[colorlinks,linkcolor=black,citecolor=blue,urlcolor=blue]{hyperref}
\usepackage[dvips]{graphics}
\usepackage{graphicx}
\usepackage{subcaption}
\usepackage{fontenc}
\usepackage{pifont}

\usepackage[utf8]{inputenc}
\usepackage[colorlinks,linkcolor=black,citecolor=blue,urlcolor=blue]{hyperref}
\usepackage{graphicx}
\usepackage{todonotes}
\usepackage{authblk}

\usepackage{appendix}
\usepackage{chngcntr}
\usepackage{etoolbox}
\usepackage{lipsum}

\newtheorem{theorem}{Theorem}

\newtheorem{lemma}{Lemma}

\theoremstyle{definition}
\newtheorem{definition}{Definition}
\newtheorem{assumption}{Assumption}

\DeclareMathOperator*{\argmax}{argmax}

\usepackage{thmtools,thm-restate}

\newtheorem{rem}{Remark}

\usepackage{todonotes}

\usepackage{bm}

\newcommand{\blue}[1]{\textcolor{blue}{#1}}

\newcommand{\norm}[1]{\left\lVert#1\right\rVert}

\newcommand{\Nats}{\mathbb{N}}

\newcommand{\R}{\mathbb{R}}

\newcommand{\E}{\mathbb{E}}
\newcommand{\Z}{\mathbb{Z}}

\newcommand{\mb}[1]{\mathbb{#1}}
\newcommand{\brac}[1]{\left(#1\right)}
\newcommand{\cbrac}[1]{\left\{#1\right\}}
\newcommand{\sbrac}[1]{\left[#1\right]}

\newcommand{\TV}{{\scriptscriptstyle T \mspace{-1mu} V}}
\newcommand{\mix}{{\scriptscriptstyle m \mspace{-1mu} i \mspace{-1mu} x}}
\newcommand{\MW}{{\scriptscriptstyle M \mspace{-1mu} W}}
\newcommand{\floor}[1]{\lfloor #1 \rfloor}
\newcommand{\mf}[1]{\mathbf{#1}}

\usepackage{todonotes}

\newcommand{\Beginproof}[1]{\begin{proof}{(Proof{#1}).}}
\newcommand{\Endproof}{\end{proof}}

\title{Markov Decision Processing Networks}
\author[1]{Sanidhay Bhambay}
\author[1]{Thirupathaiah Vasantam} 
\author[1]{Neil Walton}
\affil[1]{Durham University}

\begin{document}
\maketitle

\abstract{
We introduce Markov Decision Processing Networks (MDPNs) as a multiclass queueing network model where service is a controlled, finite-state Markov process. The model exhibits a decision-dependent service process where actions taken influence future service availability. Viewed as a two-sided queueing model, this captures settings such as assemble-to-order systems, ride-hailing platforms, cross-skilled call centers, and quantum switches.

We first characterize the capacity region of MDPNs. Unlike classical switched networks, the MDPN capacity region depends on the long-run mix of service states induced by the control of the underlying service process. We show, via a counterexample, that MaxWeight is not throughput-optimal in this class, demonstrating the distinction between MDPNs and classical queueing models.

To bridge this gap, we design a weighted average reward policy, a multiobjective MDP that leverages a two-timescale separation at the fluid scale. We prove  throughput-optimality of the resulting policy. The techniques yield a clear capacity region description and apply to a broad family of two-sided matching systems.
}%

\section{Introduction}
We introduce Markov Decision Processing Networks (MDPNs), a class of multiclass queueing networks with decision-dependent service availability.
In contrast with
Stochastic Processing Networks: current service decisions alter future service options. We characterize the capacity of these multiclass queueing networks and develop scheduling policies that maximize throughput by making decisions that better enable future service.

To illustrate the decision-dependent service that we have in mind, consider the following examples. An assemble-to-order inventory system selects inventory that is combined to fulfill an order. These service decisions reduce inventory levels, which in turn affect future availability. Thus, the optimal choices depend not only on the waiting requests but also on the evolution of inventory levels and the re-ordering process.
Ride-hailing platforms provide a further example. When the platform assigns a taxi to a customer the taxi becomes unavailable, and when it becomes free again, it will typically emerge in another part of the system. Thus, scheduling decisions (of taxis to ride-hail requests) directly impacts the network's ability to meet future demand and, therefore, has implications on future decision-making.
Many call centers implement skill based routing, here a controller must determine which staff members are assigned to handle each type of call. Draining the pool of skilled staff may limit our ability to effectively serve future calls. Thus, again, service decisions impact future customers.
The original motivation for this work is in quantum communications. In a quantum communication network, entangled photons are communicated across the network. The measurement of these particles transports classical information but also destroys the quantum information between them. An implication is that switches and repeaters that enable communication, will need to generate and storage of new entanglements to serve incoming communication requests. The controller must decide how to store and measure qubits which can decohere over time. This creates a similar trade-off between immediate and future communication capabilities.

In each of the cases above, future service is depends on current decisions. By contrast, for scheduling in traditional multiclass queueing network models, a server typically refreshes upon job completion. However, this is not the case in the examples above, a key observation is that service availability forms a Markov decision process (MDP). Although, as we will see, for throughput optimization, we must consider it as a multi-objective MDP. Consequently, optimal scheduling algorithms require ideas from both stochastic processing networks and Markov decision processes. For this reason, we refer to this queueing network model as Markov Decision Processing Network (MDPN).

In this article, we characterize the capacity region of MDPNs; we show via a counterexample that a standard scheduling policy, MaxWeight, is not throughput-optimal for MDPNs, and we construct a throughput-optimal policy derived from a weighted average-reward MDP policy (WARP). We now provide a more detailed description of our model and results. 


\subsection{Markov Decision Processing Networks}

We now give an informal description of a Markov Decision Processing Network. (A formal definition is given in Section \ref{sec:system_model}) We also discuss how this queueing model naturally describes two-sided matching networks.
We provide a comparison between our MDPN setting and switched queueing networks. In particular, we discuss the implications for the stability region of these systems. We then develop the stability and throughput properties of these queueing networks by viewing their service in terms of vector-valued MDP.

\subsubsection{Markov Decision Processing Networks and  Two-sided Matching Models.}
\begin{figure}[h!]
    \centering    
    \begin{subfigure}[b]{0.48\textwidth}
        \centering
        \includegraphics[width=\textwidth]{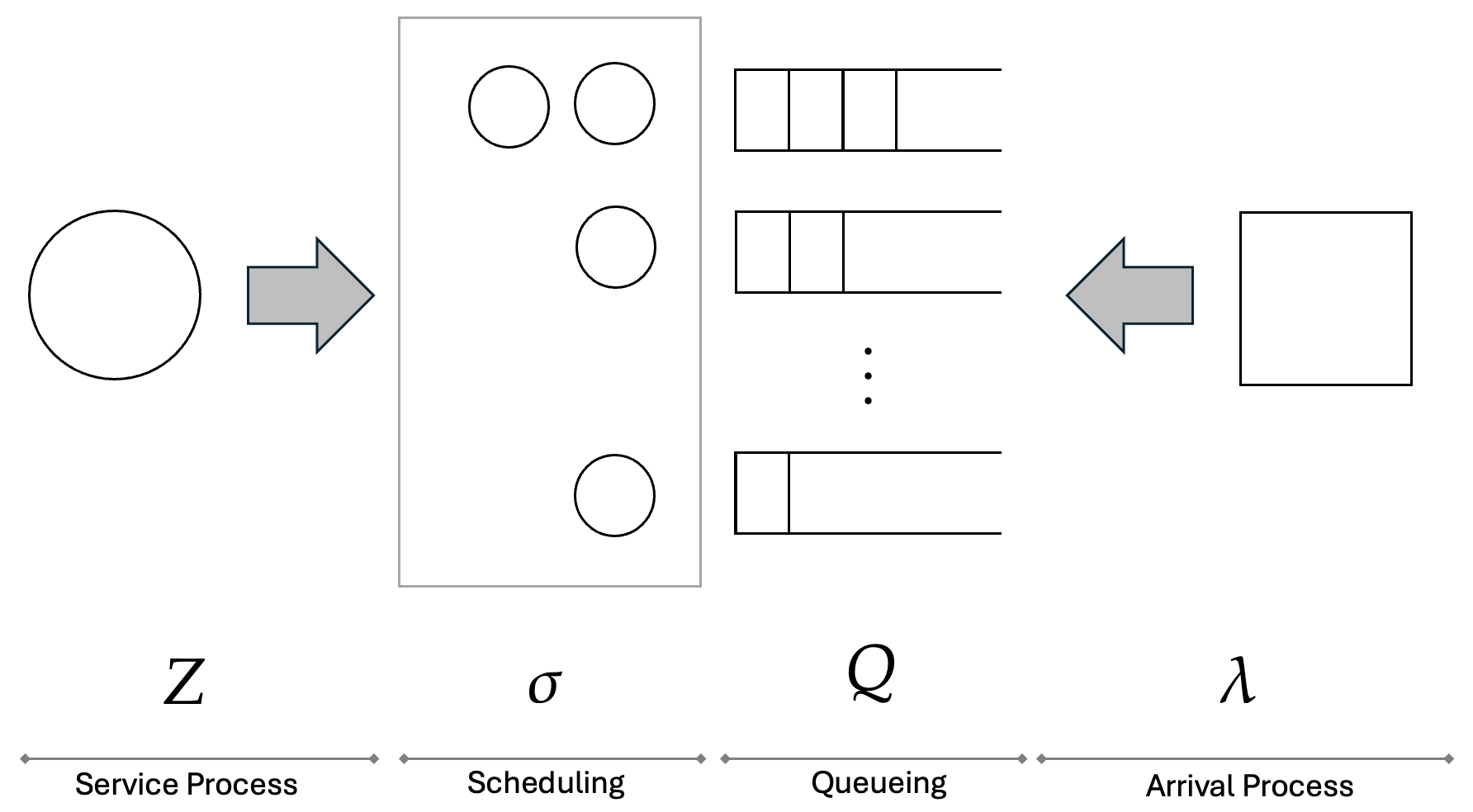}
    \end{subfigure}
    \hfill
    \begin{subfigure}[b]{0.48\textwidth}
        \centering
        \includegraphics[width=\textwidth]{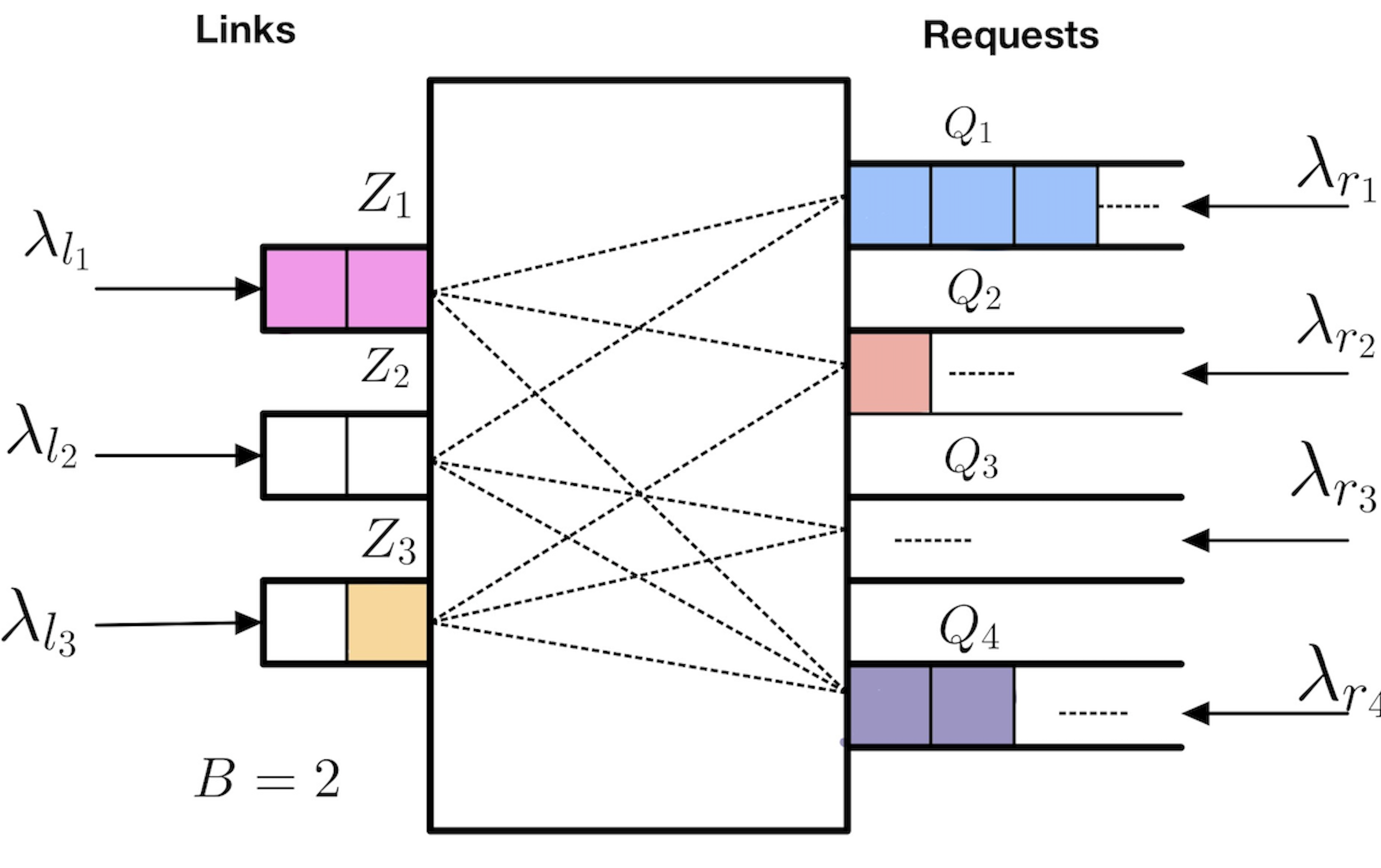}
    \end{subfigure}
    \caption{(a) An abstract representation of a Markov Decision Processing Network. The service process is an endogenous finite-state MDP, whose actions determine scheduling decisions, which in turn impact the queueing process. Arrivals are exogenous. The arrival and service processes determine the queueing process.
 (b) \label{Fig:compy} A specific example of a Quantum Switch. Entanglements are generated across links and stored in memory. Entanglements decohere and can fail when attempting to serve requests. The entanglements serve communication requests. Schedules are a matching between a request and a set of link-level entanglements.
    }
\end{figure}

An MDPN is a discrete-time multi-class queueing network model, which we represent in Figure \ref{Fig:compy}(a). 
Here, there are two key components: its service process and its arrival process, which then determine the queueing process. The service process is a finite state MDP. The actions of this MDP are scheduling decisions. The state and scheduling decisions determine the next state of the service process. These MDP dynamics represent the endogenous dynamics of the service process. The arrival process is an exogenous process that generates jobs in the queueing process. Jobs of each class queue in an infinite capacity FIFO buffer until served. The number of arrivals and the number scheduled by the MDP then determine the change in the queueing process.

An MDPN can often be seen as a two-sided matching queue. See Figure~\ref{Fig:compy}(b) for a specific example, which represents a switch in a (quantum) communication network. 
Here, categories of servers are queued in links [on the left], and jobs join queues for each request type [on the right]. The number of requests that can be queued is unbounded, while the number of services the system can support is bounded. The scheduling of a set of requests not only removes jobs from the right-hand queues but also induces a transition in the set of available servers on the left-hand side. 
Here we see that when the state of the MDP is a queueing process, then the resulting network is a queueing system with two-sided matching.

\subsubsection{Comparison with Switched Queueing Networks.}
We briefly contrast the MDPN model with switched queueing networks.

Switched queueing networks are a discrete-time model of server resource contention (see \cite{georgiadis2006resource} and the literature review below). At every time step, a scheduling decision is made to serve jobs across different queues. The set of feasible schedules is fixed and available to the scheduler at every time step.

The stability region of a switched queueing network is the set of arrival rates for which the network is positive recurrent. For switched queueing networks, the stability region is well-known to be the convex combination of the feasible schedules. 
A policy is said to be throughput-optimal if it is positive recurrent for all arrival rates within this region.
The MaxWeight policies, introduced by \citep{tassiulas1990stability}, are known to be throughput-optimal, stabilizing the network without prior knowledge or explicit estimation of arrival rates.

MDPNs exhibit fundamental differences. Here the feasible service states are generated by a controlled, Markovian service process: each scheduling decision both serves work now and changes the service state, thereby altering the options available in the future. As a result, the ability to stabilize a given arrival rate depends on the long-run mix of service states induced by the control policy, not just on the instantaneous schedules available at any given time. This makes capacity characterization more subtle for MPDNs.

Further these differences have algorithmic consequences. A direct application of MaxWeight does not lead to throughput-optimality in MDPNs; we provide a counterexample in Theorem~\ref{thrm:MWInstable}. We later recover throughput-optimality by viewing the control of the service process as an adaptive optimization process and designing a weighted, average‑reward control rule to achieve throughput optimality.

\subsubsection{Vector-valued Markov Decision Processes and Stability.}
One might notice that we have not defined a reward function for the MDP representing the service process. 
This is deliberate; we do not summarise overall system behaviour in terms of a single reward, a point that we will develop further in the paper.
If we look at the system's service process, ignoring any interactions with the request queues, then the service the request queue receives depends on its current state and past service decisions. Such a process is a Markov decision process. However, the MDP is non-standard because service impacts multiple job classes, and the service received is a vector. In effect, a single Markov decision process provides a vector of rewards for each request type.  
So a scheduling decision induces transitions in multiple queues; thus, we view the service process as a vector-valued MDP with a reward for each queue in our multiclass queueing network. 
Instead of maximizing a single reward function, we ask how we can construct scheduling rules that result in desirable queueing behavior for this vector-valued reward process. 
This observation that the service process is a vector-valued MDP is critical to understanding the stability of these systems. 


The above interpretation, then, leads to a natural formulation of the stability region of an MDPN.
An MDPN is stable if a policy has a stationary reward vector that strictly dominates the request arrival rate. Conversely, if no such policy exists, the queueing network will be unstable. The formal proof of this is one of the main results of this paper. (See Theorem \ref{thrm:stabilityRegion}). With the above observation that the service process is a vector-MDP, the stability characterization is intuitively clear. However, as we discuss in the literature review below, without the MDP interpretation, the stability analysis of general two-sided switch models was considered to be a significant challenge in prior works, with stability proofs limited to specific network topologies.

\subsubsection{A Throughput Optimality Policy for Markov Decision Processing Networks.}

If we consider the dynamics of an MDPN, we could imagine that a controller would be to estimate and optimize indigenous dynamics, such as the service process. However, the arrival process is exogenously defined and thus would not naturally be accurately estimated under an adaptive control scheme.

The joint coupling of the arrival and service processes in an MDPN suggests that we need to know the rate at which jobs of each type arrive to construct a policy that stabilizes the network. Also, we find standard arrival-agnostic policies, such as MaxWeight policies, are not throughput optimal in our setting. Given the MDP interpretation above, we need to know all the state-transition dynamics, including the arrival process. An optimal MDP typically requires knowledge of all transition probabilities. However, this is not the case here; throughput-optimality can be achieved by monitoring the queue lengths of request queues rather than estimating arrival rates.

The intuition for this is as follows. The number of jobs in the request queues can grow indefinitely while the pool of servers remains finite. Thus, when appropriately scaled, we find that many transitions can occur in the state of the server system to produce a small relative change in the request queues. In other words, there is a natural separation between the timescale of the request queues and the evolution of the service process. Thus, when busy, the request queues in effect only see the stationary dynamics of the service process. This decouples the two processes: the resulting timescale separation enables us to analyze policies that, over a given time window, are insensitive to changes in the number of jobs in the request queues and the arrival process.

The consequence is that to have a throughput-optimal scheme, a policy must solve an average-reward MDP, whose reward function is the weighted sum of queue sizes. The rationale behind this optimization is a Lyapunov drift argument, similar to the MaxWeight optimization; however, we must optimize a MDP rather than a myopic set of schedules.
In effect, we must optimize over the MDP policy space rather than over a discrete set of (instantaneous) schedules.
This argument gives the required throughput optimality: although arrival rates may vary, we do not need to estimate them to ensure stability; we only need to understand the system's internal service dynamics. The corresponding result is proven in Theorem \ref{thm:stochastic_stability}.

\subsection{Contributions}

We outline key contributions in this work:
\medskip
\begin{enumerate}
    \item To better model control in modern stochastic server systems, we define a queueing model called Markov Decision Processing Networks (MDPNs),  which generalize switched queueing networks.
    \item We characterized the stability region for MDPNs, and provided counterexamples for the MaxWeight policy, the benchmark policy of multiclass queueing networks. For this reason new policies are needed to stability these networks.
    \item We propose a novel scheduling policy, the \textit{Weighted Average Reward Policy} (WARP), which is designed to achieve throughput optimality by balancing the needs of different request classes based on their arrival rates and service requirements.
    \item To establish the optimality of the WARP policy, several technical contributions are made. As the switch becomes congested, the server states will evolve on a faster timescale compared to the incoming entanglement requests. We develope a timescale separation analysis in this work to establish a fluid limit. 
\end{enumerate}
\medskip

\noindent At a methodological level, the interplay between timescale separation and MDP in optimal scheduling of matching queues is intriguing and potentially of interest to the performance evaluation of general matching systems. Our policy demonstrates novel methods for characterizing and optimizing capacity. These results resolve several problems that were previously open on quantum switches \cite{Nain_switch,Wenhan_switch,zubeldia2022matching,dai2022queueing}.





\subsection{Literature Review}

We provide background on different queueing networks and MDP models as well as applications that motivate the MDPN model class. 

\subsubsection{Switch Networks.}

As discussed above, switched queueing networks are a discrete-time model with server contention. 
Unlike our setting, traditionally, servers are always available for the next job in a typical queueing model. In this setting, the MaxWeight algorithm is considered one of the best policies. The MaxWeight algorithm was first introduced in the context of wireless ad-hoc networks \cite{tassiulas1990stability} and then developed for Internet Switch design by \cite{mckeown1999achieving}. \cite{tassiulas1990stability} establishes the throughput optimality of MaxWeight with a Quadratic Lyapunov. Subsequent works extending these findings, see \cite{andrews2004scheduling,shah2012switched,DaiLin08}.
Several variants of MaxWeight exist, specifically BackPressure (defined by \cite{tassiulas1990stability}), which applies to switch networks with joint scheduling and routing, and MaxPressure (defined by \cite{DaiLin08} and discussed below), which applies to stochastic processing networks with continuous time and IID random service. In this paper, we focus on the scheduling problem. We will henceforth use the term MaxWeight.

When we consider models outside the framework of switched networks, the throughput optimality of MaxWeight may be lost, for instance, see \cite{van2009instability,bramson2021stability}.
We demonstrate, via a counterexample, that it exhibits suboptimal throughput behavior in settings with decision-dependent service availability. Such counterexamples highlight the non-trivial stability behavior of the MDPN model.

\subsubsection{Stochastic Processing Networks.}
Stochastic processing networks,  introduced by \cite{harrison1999heavy} and \cite{harrison2003brownian}, are a generalized class of queueing networks that account for multiple activities performed on different queues. Stochastic Processing Networks allow for diverse interactions between jobs and resources, making them suitable for a wide range of applications.

Versions of the MaxWeight and Backpressure policies have been created to stabilize SPNs, such as the MaxPressure policies, see \cite{dai2005maximum,DaiLin08,ata2008heavy}. 
The MaxPressure policy is throughput optimal in certain conditions, providing a robust framework for managing complex queueing networks. Nonetheless, as discussed above, it remains unstable in the MDPN setting.

The Markov Decision Processing Networks generalize the notion of service considered in a Stochastic Processing Network by allowing for Markov Decision dependent service activity. However, we assume that jobs leave the system after service, and so our routing structure is potentially simpler than that considered in a SPN model.
Extending this work to SPNs would be an interesting direction for future research. We refer the reader to \cite{dai2020processing}, a broad presentation of Stochastic Processing Networks, and for further discussion on the MaxPressure policies.

Although we focus on MDP formulations the connections with reinforcement learning are clear. Here we could imagine allowing the MDP time to be learned in addition to the time that MDP requires to mix thus leading to the joint learning and optimization of throughput in the MDP. Recent works focusing on optimization and learning in queueing networks include \cite{liu2019reinforcement,chen2023online,dai2022queueing}

\subsubsection{Vector Markov Decision Processes and Multi-Objective MDPs.}
In this work, we draw connections with vector-valued MDPs.
Vector-Valued Markov Decision Problems are often referred to as Multi-Objective Markov Decision Processes. These have a long history going back to the work of \cite{WHITE1982639} and \cite{furukawa1980characterization}.  \cite{roijers2013survey} provides a comprehensive survey of these methods and algorithms.  Also see \cite{hayes2022practical} for a survey of recent techniques and applications in reinforcement learning.

In terms of queueing applications, the book of \cite{altman2021constrained} considers a multi-objective approach via constrained MDPs; there are MDPs with a single reward function but constraints on multiple. 
While there is extensive work on MOMDPs, there are few works that directly investigate methodological connections between multi-objective MDPs and queueing theory.

\subsubsection{Example Applications.}
We review queueing systems that have random service availability and that could be modelled with an MDPN. As discussed Quantum Communications provided the original motivation for this work and further examples are intended to be indicative of further applications that might benefit from this framework.

\textit{Server Pools and Two-sided Matching Queues.}
Frequently, queueing systems are considered as settings where there are customers and server pools that must be matched. For instance, see \cite{bambos1993scheduling}. Specific structures are considered and a reviewed in the context of call centres by \cite{gans2003telephone}.

\cite{caldentey2009fcfs} and subsequent works such as \cite{adan2012exact} and \cite{adan2018reversibility} consider the impact of instantaneous matching in two-sided queueing systems.
Algorithms for controlling matching queues are presented in \cite{gurvich2015dynamic} and \cite{mandelbaum2004scheduling}.

\textit{Assemble-to-Order Networks.}
In an assemble-to-order system, partially assembled components are kept in inventory and combined when orders are requested. This then creates a coupled control problem with the replenishment and allocation of orders. The features of these systems are surveyed in by \cite{song2003supply}.

In assemble-to-order systems, current orders and allocations reduce future component inventories, coupling immediate fulfillment with replenishment lead times. The work of \cite{harrison1973assembly} considers assemble-to-order systems as two-sided queueing systems, emphasizing the instability that can occur and the need for non-trivial interaction between inventory (service) and orders (customers).

\textit{Ride-Hailing.}
The connections between matching, queueing, and taxis go back at least to \cite{kendall1951some}, and have attracted considerable interest with the rise of mobile platforms enabling flexible demand and supply management. Recent works that focus on two-sided matching include \cite{aouad2020dynamic} (time-limited matching decisions), \cite{besbes2022spatial} (server positioning), and \cite{banerjee2018state} (MaxWeight variants to manage demand imbalance).

\textit{Skilled Routing in Call Centres.}
Call centres with multi-skilled agents provide a further example of a system where service availability is decision-dependent: how calls are routed today shapes which types of servers are available in the future. Incoming calls belong to various categories, and each agent may be trained in a subset of those categories. \cite{wallace2005staffing} shows that maintaining a pool of cross-trained staff significantly improves performance.

An early survey of the framework and models is \cite{gans2003telephone}.
\cite{chen2020survey} provides a more recent study with a particular focus on healthcare. They also offer a series of MDP interpretations in line with the MSPN framework presented here.
Connections between skill-based routing and decision processes are emphasized in \cite{Roubos_Bhulai_2012}. And the recent work of \cite{van2024learning} considers the implications of reinforcement learning in skill-based routing choices.

\textit{Quantum Networking.}
In a quantum network, entanglement is created between users or between users and repeaters/switches for longer-distance communication. These entanglements are called link-level entanglements (LLEs). LLEs can be stored or measured (and thus destroyed). These also decay over time. This creates a two-sided matching problem between LLEs and requests for communication. For background, see \cite{HK_LO_qinternet,calefi_qi}.

There is considerable interest in quantum switches across computer science \cite{zubeldia2022matching,Longbo}, electrical engineering \cite{valls2024brief,promponas2023full}, physics \cite{collins2005quantum}, and industry\footnote{See also: \url{https://www.psiquantum.com/blueprint} and \url{https://www.youtube.com/watch?v=U5pRnK7dGcI}} \cite{bartolucci2021switch,mandil2023quantum,palacios2018introduction}. However, the queueing-theoretic view of a quantum switch as a two-sided stochastic network is still developing; our work adopts precisely this perspective and analyzes the resulting decision-dependent service dynamics.

\subsection{Organization}
In Section~\ref{sec:system_model}, we define our primary model. This section introduces basic mathematical notation, the Markov Decision Processing Network model, and a description of its dynamics.
In Section~\ref{sec:policies}, we discuss scheduling policies, demonstrate that MaxWeight is suboptimal, and define the capacity region of these Markov chains as a function of the stationary server process. We further discuss timescale separation and use this to give an optimal policy for these systems. Our theoretical findings are presented in Section~\ref{sec:main_results}. These findings consist of a characterization of the fluid limit of a quantum switch. The proof of this limit requires a timescale separation analysis. We demonstrate the optimality of the fluid limit of our policy, and we prove the throughput optimality of our policy.

\section{Markov Decision Processing Network: Model description}


\label{sec:system_model}
This section provides a formal description of the model and dynamics of an MDPN. We first introduce basic notation, then describe the network topology, and finally present the dynamics of the system.

\subsection{Basic Notation}
We denote  the set of non-negative integers as $\Z_+$ and  $[M]$ corresponds to $\{1,\ldots, M\}$. We use ${\mathbb R}$ to denote the set of extended reals, i.e., $\mathbb R\cup \{\infty\}$. For $x,y \in \mathbb R$, $x\vee y = \max\{x,y\}$ and $x\wedge y=\min\{x,y \}$. For $i \in [M]$, $\bm e_i$ is the $i$-th unit vector in $\mathbb R^M$; $\bm 1$ is the vector of all ones; $\bm 0$ is the zero vector. Given vectors $\bm u, \bm v \in \mathbb R^M $, we write $\bm u \leq \bm v$ if $u_i \leq v_i$ for $i\in[M]$. For a set $\mathcal{M}$, we use $|\mathcal{M}|$ to denote its cardinality. We use $\bm a \cdot \bm b$ to denote the dot product and $\bm a \bm b$ to represent the element-wise multiplication between two vectors of the same dimension. We use the notation $[y]_+=\max(y,0)$ and $[y]_+^a=\max(y,a)$.

\subsection{Network Model}
\label{subsec:topology}
Here we describe the queueing process, then the service process and finally the joint Markovian dynamics of an MDPN.



\subsubsection{Requests: Queues and Arrivals}
We consider a discrete-time queueing model where, at the beginning of each time slot, requests of $R$ types arrive. Let $\mathcal{R}=\{1,\dots,R\}$ denote the set of request types. Requests of type $r\in\mathcal{R}$ are stored in an infinite-capacity FIFO queue.

Denote by $Q_r(t)\in\mathbb{Z}_+$ the number of type-$r$ requests present at the beginning of time slot $t$ (just before scheduling). We define the queue-length vector
$\bm Q(t)=(Q_r(t):r\in\mathcal{R})$. When $\sigma_r(t)$ requests are to be served for class $r$ and $\lambda_r(t)$ new arrivals occur in the slot, the queue dynamics satisfy
\begin{equation}
Q_r(t+1) = \big[\,Q_r(t)-\sigma_r(t)\,\big]_+ + \lambda_r(t),\qquad r\in\mathcal R .
\label{eq:queue_update}
\end{equation}
Type $r$ request arrive according to an exogenous arrival process $A_r(\cdot)$ whose increments $\lambda_r(t) := A_r(t)-A_r(t-1)$ are i.i.d. with mean $\lambda_r \in \R_+$ and are positive and bounded. 

\subsubsection{Servers: States, Actions and Schedules.}
The (finite) set $\mathcal Z$ represents the server states. At each decision epoch, the server process $\bm Z$ evolves according to a Markov decision process. Specifically, given the current server state $\bm z \in \mathcal Z$ and a control action $a$, the next server state $\bm z'$ is determined by the probability 
\[
P(\bm z' | \bm z, a) = \mathbb{P}\left(\bm Z(t+1) = \bm z' \mid \bm Z(t) = \bm z, a(t) = a\right).
\]
We let $\mathcal{A}$ be the (finite) set of actions.

Note that the above describes the transitions of an MDP; however, we do not currently give a reward function. In the context of MDPNs, the reward function is will be defined to reflect the number of requests served for each class, and is used to optimize throughput of the request queues.

We now extend the above transitions so that we can schedule the request queues.
A schedule $\bm \sigma = (\sigma_r : r \in \mathcal R) \in \mathbb Z^R_+$ specifies the number of requests scheduled for service in a given time slot, based on the current server state $\bm z$ and the chosen action $a$. 
We let $\mathcal{S}$ denote the set of all feasible schedules.
We make the following monotonicity assumption: if ${ \sigma} \in \mathcal S$ can be realized as the average schedule under some action, then so can every non-negative vector with $\tilde {\sigma} \leq {\sigma}$.
  Or, informally stated, it is always possible to do less work on any given set of request queues.
The joint transition probability for the next server state, $\bm z'$ and the selected schedule $\bm \sigma$ is given by
\begin{align}\label{eq:joint_transition}
P(\bm z', \bm \sigma \mid a, \bm z) &= \mathbb{P}\left(\bm Z(t+1) = \bm z',\, \bm \sigma(t) = \bm \sigma \mid a(t) = a,\, \bm Z(t) = \bm z\right) \\
P(\bm \sigma \mid a, \bm z) &= \sum_{\bm z'} P(\bm z', \bm \sigma \mid a, \bm z) \, .
\end{align}
We also let ${\bm \sigma}(a, \bm z)$ denote the average schedule under action $a$ in state $z$, i.e.,
\[
{\bm \sigma}(a, \bm z) = \mathbb{E}[\bm{\sigma} \mid a, \bm z] = \sum_{\bm{z}' \in \mathcal{Z}} \sum_{\bm{\sigma} \in \mathcal{S}} \bm{\sigma} \, P(\bm{z}', \bm{\sigma} \mid a, \bm{z}).
\]


\subsubsection{Order of events.}

Over one time slot,
the following sequence of events occurs. First, the server state becomes available.
Then, an action is chosen (typically, an action will be selected as a function of the server state and the request queue lengths).
The server state and action induce a transition that produces both the next server state and the schedule that can be implemented at the request queues, cf. \eqref{eq:joint_transition}.
The request queues then evolve according to the selected schedule, as given in Equation \eqref{eq:queue_update}.

\subsubsection{Markov Process Dynamics.}
We define the joint process $\bm X(t) := (\bm Q(t),\bm Z(t))$
with
$\bm X =(\bm X(t):t\in \Z_+)$. We let $\mathcal X$ denote the set of states of this process.
For the policies considered in this paper, $\bm X$ will be a Markov process.
 In particular, we will assume that the action $a(t)$ chosen at time $t$ is a function of the current state $\bm X(t)$. (I.e. restrict ourselves to stationary policies as this tradition in the analysis of MDPs.)
When the system is positive recurrent, we use $\bm X(\infty) = \brac{\mf Q(\infty),\mf Z(\infty)}$
to denote the steady-state values. 

\subsection{Policies}
A policy, which we will denote by $\pi$, is a rule that determines the action $a$ based on the past observed state of the system.

We refer to a policy as Markovian if it depends only on the current state $\bm X(t)$. We also allow for randomized choice schedules. With this in mind, we define $\langle \mathcal A \rangle$ to be the set of random variables with support on $\mathcal A$. Thus, a policy is a function from the set of states to the set of (randomized) schedules $\bm \pi : \mathcal X \rightarrow  \langle\mathcal A \rangle $. 
For $\bm x \in \mathcal X$, we will apply the notation:
$\pi(a \mid \bm{x})$
To denote the probability that action $a$ is selected when the MDPN is in state $\bm x$ under policy $\pi$.

We say a policy is \emph{request-agnostic} or \emph{agnostic} if it does not have knowledge of the number of request queue lengths. Specifically, a request-agnostic policy is a function from the set of server states to the schedules, $\bm \pi : \mathcal Z \rightarrow \langle\mathcal A \rangle$.
We let $\mathcal P$ denote the set of request-agnostic policies.
Similarly, we define $\bm \pi(a |\bm z) $ to be the probability that action $a$ is selected when the MDPN is in state $\bm z$ under policy $\bm \pi$.

\section{An Optimal Scheduling Scheme}
\label{sec:policies}

For MDPNs, one of the objectives of scheduling policies is to maximize throughput—the rate at which requests are successfully served. 
However, the scheduling policies considered in this work must not only maximize the throughput of request types but also effectively manage the allocation of servers. This approach is distinct from prior work on scheduling in switched queueing networks. Notice
first, if we were to solve an MDP directly, we would typically need to know arrival rate parameters $\bm \lambda$. Thus, the policy will not typically be throughput optimal, as we would require a policy for every arrival rate vector. Second, as will prove shortly,  directly applying the standard throughput optimal policy leads to instability in MDPs. 

However, before proceeding further, we note that it is not even clear what throughputs are achievable for an MDPN; we characterize this in the following section.

\subsection{Capacity Region}
The capacity region is the set of request rates for which the system can be stabilized. More formally, the capacity region, $\mathcal{C} $, is the set of arrival rates for requests $\bm \lambda \in \mathbb R _ {+} ^R$ for which there exists a policy where the queue size process $\bm X(t)$ is a positive recurrent Markov chain.

When server availability is independent of the system state, there are standard arguments that characterize the capacity region, as seen in \cite{andrews2004scheduling}. However, this is not the case for our matching system; there is an interdependence between both sides of the queueing system. Thus, the characterization of the capacity region below is new and non-standard because service on the right-hand queues is interdependent with the evolution of the left-hand-side queues. The notation of a request-agnostic policy is key to characterizing the capacity region.

\begin{restatable}{theorem}{ThrmStabilityRegion}\label{thrm:stabilityRegion}
 A necessary condition for $\bm \lambda \in \mathcal C$ is that there exists a request agnostic policy $\bm \pi $ with server process having stationary distribution $\mu$ such that 
\begin{equation}\label{eq:Ccirc}
\lambda_r <  \sum_{\bm z \in \mathcal Z} \sum_{a\in\mathcal{A}}  \mu(\bm z)\pi(a | \bm z) \sigma_r(\bm z, a)  \, , \qquad \forall r  \in \mathcal R \, .   
\end{equation}  
A sufficient condition for $\bm \lambda \in \mathcal C$ is that there exists a request agnostic policy $\bm \pi$ with server process having stationary distribution $\mu$ such that 
\begin{equation}\label{eq:Ccirc2}
\lambda_r \leq  \sum_{\bm z \in \mathcal Z} \sum_{a\in\mathcal{A}}  \mu(\bm z)\pi(a | \bm z) \sigma_r(\bm z, a)  \, , \qquad \forall r  \in \mathcal R \, . 
\end{equation}
\end{restatable}

The proof of Theorem \ref{thrm:stabilityRegion} is given in Section~\ref{append:stabilityregion} in the appendix.
The above result characterizes the set of arrival rates that can be stabilized. 
We define $\mathcal {C}^\circ$ to be the set of arrival rates such that \eqref{eq:Ccirc} holds. 
We now provide the definitions of throughput optimality.
\begin{definition}
Policy $\bm \pi$ is \emph{throughput optimal} if it's positive recurrent for all arrival rates $\lambda\in \mathcal C^\circ$. 
\end{definition}

 Informally stated, a policy is throughput optimal if it is stable for the maximum set of arrival rates.
 We note that the set of stabilizable policies depends on the set of stationary distributions of the left-hand queueing system. To optimize service, we must optimize over the stationary distribution of the servers. This provides one justification for optimizing an average reward MDP for the left-hand queues to achieve throughput optimality.


\subsection{MaxWeight Scheduling}
The most widely studied policy in the network community is the MaxWeight. It selects the schedule that maximizes a weighted sum of queue lengths in each time slot.
\begin{definition}(MaxWeight) 
The MaxWeight policy $\pi^{\MW}$ selects the action that solves the following optimisation:
\begin{equation}
\label{eqn:MW_schedule}
 \max_{ a \in \mathcal A}\sum_{r\in \mathcal{R}}Q_r {\sigma}_r(a,\bm z)
\end{equation}
when the network is state $\bm z \in \mathcal Z$  with queue lengths $\bm{Q}= (Q_r : r \in \mathcal{R}) $. 
\end{definition}

For an input-queued switch (described by \cite{mckeown1999achieving}), the MaxWeight policy is known to be maximally stable in that it stabilizes the system if the average arrival rates $\bm \lambda=\brac{\lambda_r,r\in\mathcal{R}}$  lie within the capacity region. 
This observation has been extended to settings where the set of feasible schedules is time-varying as a Markov process, see \cite{georgiadis2006resource}.
In these cases, the capacity region is a convex combination of the set of schedules. However, as we have seen in Theorem \ref{thrm:stabilityRegion}, it is no longer the capacity region, $\mathcal{C}$, for an MDPN. As we now discuss, this leads to a drop in performance when directly applying MaxWeight.

\subsection{A Counter Example: MaxWeight is not Throughput Optimal.}
\label{subsec:counter_example}

We now give a simple example that demonstrates that MaxWeight is not throughput optimal.

Suppose that there are three types of servers  $s_1$, $s_2$, $s_3$, and three types of request: $r_1$ which requires one servers $s_1$ only; $r_2$ which requires a server $s_2$ only; and $r_3$ which requires all servers $s_1$, $s_2$, $s_3$. 
We suppose that the servers arrive in sequence $1,2,3,1,2,3,...$, or more formally, a server $s_1$ arrives at time $3n+1$, server $s_2$ at time $3n+2$, and server $s_3$ at time $3n+3$ for $n\in \mathbb Z_+$. 
Suppose that request arrivals are a Bernoulli process with parameters $\lambda_{r_1},\lambda_{r_2},\lambda_{r_3}$ respectively and occur just before the first of these three time slots.

Notice in this setup, the $r_1$ and $r_2$ request queues are allowed to schedule the servers in $l_1$ and $l_2$ before $r_3$. If we follow a myopic policy, such as MaxWeight, then we will schedule these requests as long as the queues for $r_1$ or $r_2$ are non-empty. However, if we schedule queues $r_1$ or $r_2$, then we cannot serve requests in queue $r_3$. In other words, $r_1$ and $r_2$ have priority over $r_3$. Because we can only serve $r_3$ requests when there are no arrivals for $r_1$ and $r_2$, this induces the following necessary condition for stability under MaxWeight for the switch just described:
\[
\lambda_{r_3} < (1-\lambda_{r_1}) (1-\lambda_{r_2})\,.
\]
However, we can plan ahead slightly, and in the last of the three time slots, we can choose which requests to serve. This results in the following sufficient condition for stability :
\[
\lambda_{r_1} + \lambda_{r_3}< 1 ,\qquad \lambda_{r_2} + \lambda_{r_3} <1, \qquad \lambda_{r_3} <1 .
\]
Note that the three constraints above correspond, respectively, to the allocation of the three servers generated on links $l_1$, $l_2$, and $l_3$. Here, we can stabilize $\lambda_{r}=0.4$ $\forall r$, whereas MaxWeight is unstable.
From the inequalities above, we see that MaxWeight's stability region is strictly smaller than the capacity region for this quantum switch. 
 We summarize these findings in the following theorem:

\begin{restatable}{theorem}{MWInstable}\label{thrm:MWInstable}
    The MaxWeight policy is not throughput optimal for MDPNs.
\end{restatable}
A simulation demonstrating the MaxWeight instability and stability of our alternative policy is given in Figure \ref{Fig:lyapunov_drift}.
Also, a further stochastic example is developed in Appendix~\ref{append:instability_maxweight}.

So, MaxWeight is not throughput optimal. We will now demonstrate that this is because MaxWeight does not plan ahead. The logic behind MaxWeight remains valuable in designing throughput-optimal policies, and it can provide maximal stability in some matching systems with stochastic service, in particular \cite{zubeldia2022matching} provides MDPN models for which MaxWeight is maximally stable. However, we urge caution when directly implementing MaxWeight in settings with decision-dependent server processes.

\subsection{Timescale Separation}
\label{sec:timescale}

In this section, we discuss timescale separation in an MDPN. We also refer the reader to \cite{zubeldia2022matching} for an MDPN setting and an analysis of MaxWeight in a Quantum switch with $\textsc{Y}$ and $\textsc{W}$ matching topologies. Our discussion here helps us gain intuition that leads to an optimal policy. The formal demonstration of the observations made here is proven across Theorems \ref{thm:fluid_convergence}, \ref{thm:fluid_stability}, and  \ref{thm:stochastic_stability}.

\begin{figure}[h!]
    \centering    
    \begin{subfigure}[b]{0.48\textwidth}
        \centering
        \includegraphics[width=\textwidth]{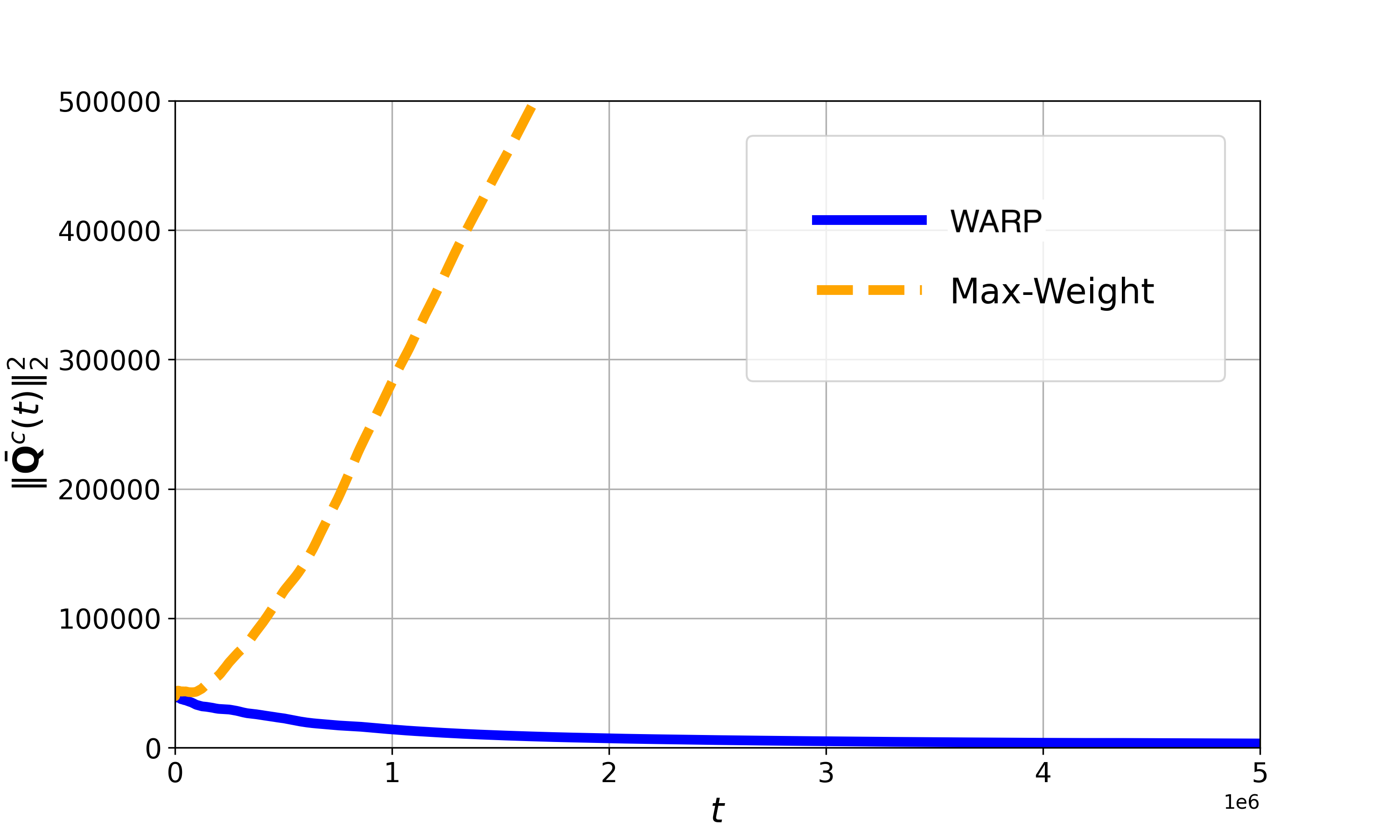}
        \caption{MaxWeight instability: We set $\lambda_{r_1}=\lambda_{r_2}=0.005$, $\lambda_{r_3}=0.004$, $\lambda_{l_1}=\lambda_{l_2}=0.02$, $\lambda_{l_3}=0.01$, $d_{l_1}=d_{l_2}=0.00001$, $d_{l_3}=0.99999$, $\gamma_{r_i}=1$ for all $i\in[3]$ and $c=200$. The optimal policy is computed through value iteration.}
        \label{Fig:lyapunov_drift}
    \end{subfigure}
    \hfill
    \begin{subfigure}[b]{0.48\textwidth}
        \centering
        \includegraphics[width=\textwidth]{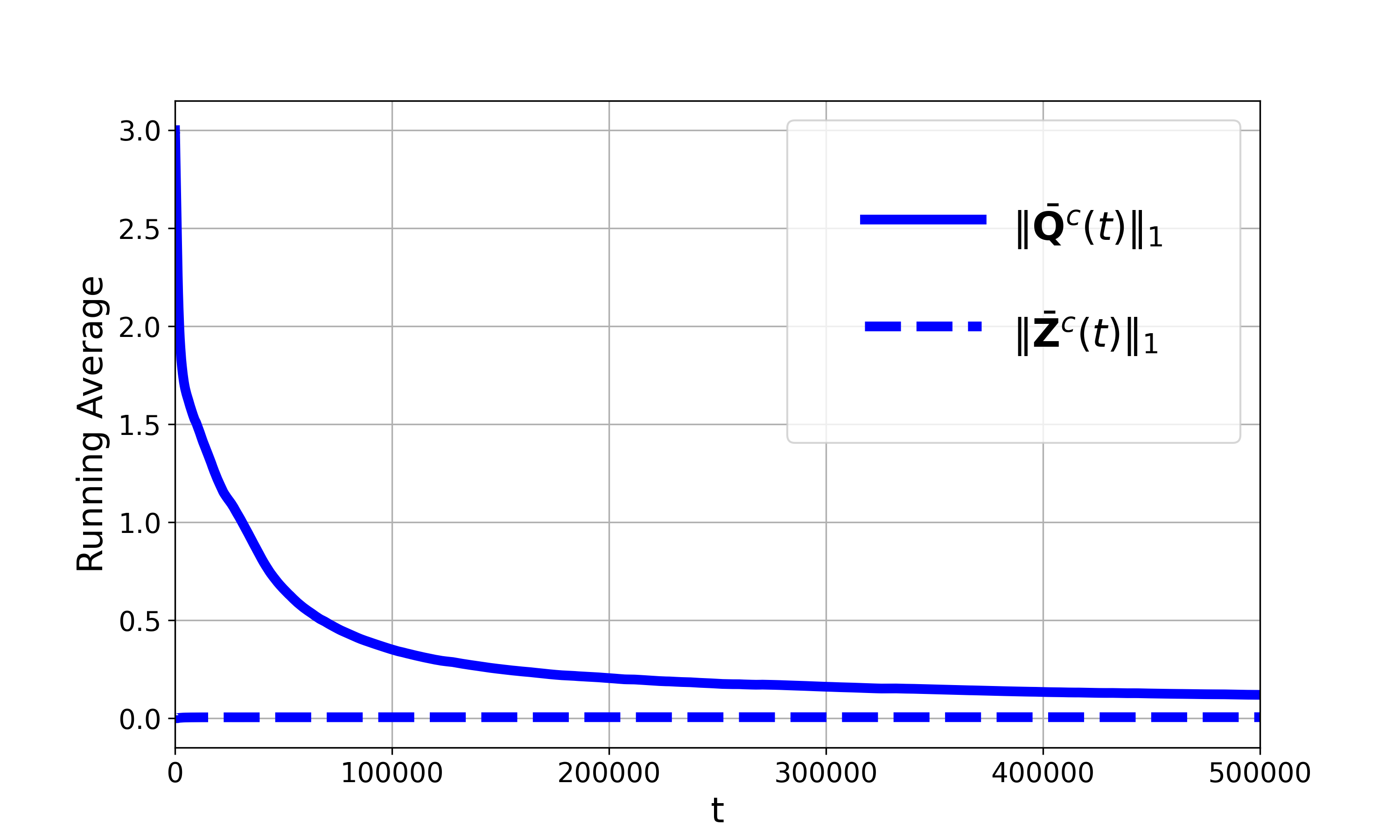}
        \caption{Timescale separation: We set $\lambda_{r_1}=\lambda_{r_2}=0.05$, $\lambda_{r_3}=0.04$, $\lambda_{l_1}=\lambda_{l_2}=0.2$, $\lambda_{l_3}=0.1$, $d_{l_1}=d_{l_2}=0.00001$, $d_{l_3}=0.99999$, $\gamma_{r_i}=1$ for all $i\in[3]$ and $c=200$. The WARP policy is used for server scheduling.}
        \label{Fig:fl-qs}
    \end{subfigure}
    \caption{For simulations, we simulated the counterexample given in Section~\ref{subsec:counter_example} with probabilistic decoherence of servers. In the counterexample, we considered three types of requests ($\mathcal{R}=\{r_1,r_2,r_3\}$) with three servers ($\mathcal{Z}=\{l_1,l_2,l_3\}$). Moreover, we have considered $\mathcal{Z}_{r_1}=\{l_1\}$, $\mathcal{Z}_{r_2}=\{l_2\}$ and $\mathcal{Z}_{r_3}=\{l_1,l_2,l_3\}$. We set $B=1$.}
    \label{Fig:comp}
\end{figure}

In Figure \ref{Fig:lyapunov_drift}, we plot the performance of MaxWeight against an alternative policy, WARP, which we will describe shortly. MaxWeight has optimal drift for the quadratic Lyapunov function in a switch network Model. However, this is not the case for an 
MDPN.
To evaluate the performance of a queueing network with Markovian servers, it is essential to understand the timescale separation that occurs under congestion.

Consider the MDPN where there are many requests, i.e., $\sum_r Q_r(0) = c$ for $c\gg 1$. In this state, the number of the servers will be far smaller $\bm Z(t) = O(1)$ because of the physical limitations of the switch having bounded storage for servers. 
Timescale separation is demonstrated in Figure \ref{Fig:fl-qs}.
To make any relative change in the request process, we require $c$ transitions of the server process. Since the server process is already close to equilibrium, being of order $O(1)$, it quickly converges to its stationary behavior. Thus, the resulting scheduling dynamics placed on the request queues are ultimately determined by the stationary behavior of the server process. (This further motivates the capacity region characterization that we already proved in Theorem \ref{thrm:stabilityRegion}.)

\subsection{Addressing the Sub-Optimality of MaxWeight}
In contrast, to more classical switched queueing networks, two-time scale separation affects the Lyapunov drift of the MaxWeight in such a way that it is no longer optimal. The original rationale of the MaxWeight policy is to achieve the maximum negative Lyapunov drift. We can informally explain this as follows:
The differential equation below approximates the request queues given in Section~\ref{sec:system_model},
\begin{equation*}
    \frac{d\bar{Q}_r(t)}{dt}
    =
    \lambda_r 
        -  
        \mathbb E_{\bm Z(\infty) \sim \mu(t) } [
            {\sigma}_r (\infty)
            ]
    , \qquad  r\in \mathcal{R}. \footnote{We more formally define $\bar Q$ and associated fluid model terms in Section \ref{sec:main_results}.}
\end{equation*}
Here the state of servers is stationary with stationary  distribution $\mu(t)$. This then induces a stationary schedule, denoted by ${\sigma}_r (\infty)$ above. If
we take the function $L(\bar{\bm{Q}}(t)) = \sum_{r \in \mathcal{R}} \bar{Q}_r^2(t)/2$, which would be the Lyapunov function typically associated with MaxWeight. Then, differentiation with the chain rule gives
\[
 \frac{dL}{dt}= \sum_{r \in \mathcal R} \bar Q_r(t)\lambda_r -    \mathbb E_{\bm Z(\infty) \sim \mu(t) } \Big[ \sum_{r \in \mathcal R}\bar Q_r(t){\sigma}_r(\infty)\Big]\, .
\]
(The above notation indicates that the server process is stationary while the queue size process is not.)
The MaxWeight policy should maximize the negative drift of the Lyapunov function. However, the MaxWeight policy does not achieve the maximum negative drift for the MDPN from the above equation. Specifically, the maximum negative drift solves the optimization  
\begin{equation}\label{eq:avMDP}
    \max_{\pi \in \mathcal P }\;\; \mathbb E_{\bm Z(\infty) \sim \mu^\pi } \Big[ \sum_{r \in \mathcal R}\bar Q_r(t) {\sigma}_r(\infty)\Big] \, .
\end{equation}
In the display above, it is evident that there is a separation in timescales between the queue size process, which depends on the current time, and the service state, which is stationary with respect to the stationary distribution that depends on the current time.
In a classical switch or a wireless network, the above optimization is a combinatorial optimization problem, such as a bipartite matching problem. However, as shown above, we must jointly optimize over schedules and their induced stationary distributions. In particular, the critical observation is that the above optimization is an MDP. We discuss this point in more detail in the next subsection.

The above argument is somewhat heuristic. To make this rigorous, we must prove that the above timescale separation and fluid limit are correct. We complete this in Theorem \ref{thm:fluid_convergence}. We also verify that the solution of the MDP gives the optimal drift, which we define below. Then, we investigate how this impacts stochastic policies that implement schedulers solving this MDP. This result is given in Theorem \ref{thm:fluid_stability}. \blue{ }

\subsection{Optimal Scheduling is an Average Reward MDP}

As introduced in the previous section, the optimal policy must solve an average reward MDP. 
For illustration, see Figure~\ref{Fig:lyapunov_drift}, which shows that the drift of the policy obtained as a solution of the average reward MDP is negative compared to that of the MaxWeight policy applied to the MDPN model, which is unstable. The negative drift under the average reward MDP solution reflects the policy's ability to more effectively balance short-term scheduling decisions with the long-term availability of resources.

We now briefly discuss average reward MDPs, providing some notation for later use. For a detailed treatment of average  MDPs, we refer the reader to Chapter 7 of \cite{puterman2014markov} or Chapter 5 of \cite{bertsekas2011dynamic}.

The single-step reward associated with implementing action $a$, state $ \bm z $, and request queue length vector $\bm  q$ is
\begin{equation}
\label{eqn:reward}
u(a,\bm{z} ; \bm q) = \sum_{r \in \mathcal{R}} \sum_{\bm \sigma \in \mathcal{S}} q_r \sigma_r P(\bm \sigma \mid \bm z, a) .
\end{equation}
This reward represents the weighted average number of requests fulfilled in the current time slot. Notice in terms of the MDP framework, the system's state is $\bm z$, the action chosen is $a$, and $\bm q$ is a parameter of our objective function. 
Note that we take the queue size parameter $\bm q$ as a fixed constant, so our policy is myopic in that it does not account for the impact of decisions on the evolution of request queues.


The single-step reward associated with a (stationary) policy $ \pi $ and state $ \bm{z} $ is given by:
$
u( \pi(\bm{z}),\bm{z}) = \sum_{a \in \mathcal{A}} \pi(a \mid \bm{z}) \, u(a,\bm{z}),
$.
Notice that under any policy, $\pi$ states evolve as a positive recurrent discrete-time Markov chain with a single irreducible communicating class. Thus, our average reward MDP is \emph{unichain}, which in turn implies the existence of stationary solutions to the average reward MDP \cite[Proposition 5.2.4]{bertsekas2011dynamic}.

The average reward associated with stationary policy $ \pi $ is denoted by $ R^{\pi} $ and is defined as:
\begin{align*}
R^{\pi} & = \lim_{T \to \infty} \frac{\mathbb{E} \left[\sum_{s=0}^{T-1} u( \pi(\bm{Z}(s)),\bm{Z}(s))\right]}{T}
\\
&= \mathbb E_{\bm z\sim \mu^\pi} \Big[  \sum_{r \in \mathcal R} \sigma_r^{\pi}(\infty)  q_r \Big] \, .
\end{align*}
     Above, we note that the Markov chain $\bm Z(s)$ is ergodic, and we let $\mu^\pi$ denote the stationary distribution of $\pi$. We let the random variable $\sigma^\pi_r$ denote the stationary number of requests of type $r$ processed under $\pi$.\footnote{We assume that the policy $\pi$ will continue to process requests irrespective of whether requests are queued.} Rewriting the above expression in terms of these stationary quantities, the optimal reward and the optimal policy are denoted, respectively, by

\begin{align} \label{eq:AREMDP}
R^{\star}(\bm q) 
&
:= \max_{\pi \in \mathcal P}\,\, \mathbb E_{\bm z\sim \mu^\pi} \Big[ u(\pi(\bm z),\bm{z}; \bm q)  \Big]\, , 
\\
\pi^\star(\bm q) 
&\in \argmax_{\pi \in \mathcal P}\,\, \mathbb E_{\bm z\sim \mu^\pi} \Big[  u(\pi(\bm z),\bm{z}; \bm q) \Big]\, ,
\\
\bm \sigma^\star (\bm z; \bm q) &:=  \bm \sigma(\pi^\star(\bm z ;\bm q),\bm z) \, ,
\end{align}

for $\bm q \in \mathbb Z_+^R$, $\bm z \in \mathcal Z$ and with $u(a,\bm z; \bm q)$ as above in \eqref{eqn:reward}.
In the above expression, we think of $\bm q$ as a parameter of our optimization, rather than as a state in the MDP.
For each state $\bm{z} \in \mathcal Z$, the Bellman equation for the optimal reward $R^\star$ is given by:
\begin{align}
& R^\star(\bm q) + V^{\star}(\bm{z};\bm q) \notag
\\
&= \max_{a\in \mathcal{A}} \left\{ \sum_{r \in \mathcal{R}} \sum_{\bm \sigma \in \mathcal{S}} q_r \sigma_r P(\bm \sigma \mid \bm z, a)  +\sum_{\bm{z}' \in \mathcal{Z}}  P(\bm{z}' \mid a, \bm{z}) V^{\star}(\bm{z}';\bm q) \right\},\label{eqn:bellman}
\end{align}
where $V^{\star}(\bm{z};\bm q)$ is the relative value function, representing the relative value of state $\bm{z}$ compared to the optimal-term average reward $R^\star(\bm q)$,  $P(\bm{z}' \mid a ,\bm{z})$ is the transition probability from the current state $\bm{z}$ to the next state $\bm{z}'$, given action $a$.
The optimal policy maximizes~\eqref{eqn:bellman}, which by design balances the allocation of servers to satisfy current requests while also considering the impact on future system states. 
The optimal policy can be calculated through value iteration or policy iteration.  If the server process is not known, parameters can be learned directly or via tabular reinforcement learning algorithms, and these estimations can then be applied to calculate optimal scheduling decisions.

\subsection{Weighted Average Reward Policy (WARP)}

This paper considers the optimal policy for allocating servers in an MDPN, which is found by solving an average reward MDP that optimizes the utilization of servers to fulfill incoming requests. 

We call this policy the Weighted Average Reward Policy (WARP). It consists of a single parameter $\tau \in\mathbb N$ and the sequential implementation of the average reward MDP described above.
The WARP policy is defined as follows:

\begin{quote}
    \textbf{WARP}: Time is divided into epochs, the length of an epoch, $\tau(\bm Q)$, is a function of the queue size vector, $\bm Q$, at the beginning of the epoch. Throughout the epoch, we fix the policy to be the optional solution $\pi^\star(\bm Q)$ to the MDP \eqref{eq:AREMDP}.
\end{quote}

A key point on the WARP policy is that it does not need to know the average arrival rate $\bm \lambda$. 
The WARP policy makes decisions where current queue sizes are simply a parameter used to determine the actions of the service process regardless of the state evolution of the request queues. (I.e., if requests are requested to serve, WARP will serve them. However, even if a request queue is empty, the policy will continue to schedule jobs regardless of the availability of requests.) Since we only require current queue sizes, the WARP policy behaves similarly to the Max-Weight policy, which doesn’t require explicit knowledge of $\bm \lambda$ to achieve optimal throughput.

While traffic patterns do not need to be learned, the policy requires solving an MDP. 
MDPs are well-known to suffer from the curse of dimensionality, in our case as the buffer size grows. 
In practice, reinforcement learning methods and function approximation would be implemented offline to find the queueing network's optimal dynamics before implementation. These directions would be an important area for future investigation. See the literature review for a review of current machine learning methods for multi-objective MDPs.

\begin{assumption}\label{ass:tau}
    The function $\tau(\bm Q)$ is non-decreasing with respect to the components of the queue size vector $\bm Q$, and 
    \[
    \tau(\bm Q) = \omega\big( \log(\|\bm Q\|_1) \big).
    \]
\end{assumption}

\section{Throughput Optimality}
\label{sec:main_results}
This section presents and discusses the main results of the WARP policy, which is the following: 

\begin{theorem}
\label{thm:stochastic_stability} For any $\tau(\bm Q)$ such that Assumption \ref{ass:tau}, the WARP policy is throughput optimal.
\end{theorem}

The proof of Theorem~\ref{thm:stochastic_stability} follows from the fluid limit and stability results established in the following sections. See \cite{dai1995stability} and \cite{bramson2008stability} for an overview of this approach. However, to establish our results, we must also establish a timescale separation between the MDP process followed by servers and the states of the queues. See \cite{hunt1994large} for a related timescale separation result.

The structure of the proof Theorem \ref{thm:stochastic_stability} is as follows: 
Below in Section \ref{sec:fluid_model}, we state the fluid model associated with an MDPN and prove the throughput optimality of the fluid model. 
We then must establish the connection between the stochastic MDPN model and its fluid equivalent. 
This is stated as Theorem \ref{thm:fluid_convergence} in Section \ref{sec:fluid_model} below and proven in Section \ref{sec:fluid_limit_proof} of the appendix. 
A key aspect of establishing the fluid model and its establishment is a timescale separation between the server processes and the queue states. This is discussed in Section \ref{sec:timescale} above and stated in Proposition \ref{prop:timescale} below.
With the fluid limit and fluid model throughput optimality results in place, we can then prove Theorem \ref{thm:stochastic_stability} via the Multiplicative Foster-Lyapunov Theorem (See \cite{robert2013stochastic} and \cite{bramson2008stability}) this final part of the proof of Theorem \ref{thm:stochastic_stability} is established in Section \ref{sec:fluid_limit_proof} of the appendix.



\subsection{Fluid Model and Fluid Stability}
\label{sec:fluid_model}

We now focus on the fluid model for the WARP policy, which is defined as follows.
\begin{definition}(WARP Fluid model)
\label{defn:fluid_limit}
Given an initial condition $\bar{\bm Q}(0)\in \R_+^{R}$ such that $\norm{\bar{\bm Q}(0)}_1>0$, we say that a Lipschitz continuous function $\bar{\bm Q}:[0,\infty) \rightarrow \R_+^{R}$ is said to be a solution to the fluid model if it satisfies following equations for all $t,s\geq0$ and $r\in\mathcal{R}$ 
\begin{align}
\label{eqn:fl_1}
\bar{Q}_r(t)&=\bar{Q}_r(0)+ \bar{A}_r(t)- \bar{D}_r(t),\\
\bar{A}_r(t)&=\lambda_r t,\label{eqn:fl_2}
\end{align}
and
\begin{align}
&\int_{s}^t \sum_{r \in \mathcal R} \bar{Q}_r(u) d \bar{D}_r(u) \notag \\
&\geq \max_{{\bm \pi}\in \mathcal P} \int_s^t \sum_{r \in \mathcal R}
\sum_{\sigma \in \mathcal S}
\bar{Q}_r(u)  \sigma_r P(\bm \sigma \mid \pi(\bm z), \bm z) 
\mu^\pi(\bm z)
du
 \label{eqn:fl_3},
\end{align}
where
$\mu^{{\bm \pi}}$ denotes the steady-state distribution of the left-side process $\bm{Z}$ under policy ${\bm \pi}$ and $\bar D_r(t)$ is increasing and Lipschitz continuous.
\end{definition}

Notice that, unlike the stochastic model, which has two components $(\bm Q, \bm Z)$, in the fluid model, the state is only expressed in terms of $\bar{\bm Q}$. This is because there is a timescale separation between the two processes, and thus we can average over the stationary distribution of the server process $\bm Z$. This timescale separation is formally stated in Proposition \ref{prop:timescale}.  

\begin{definition}[Fluid Model Stability and Throughput Optimality]
A fluid model is considered to be \emph{stable} if there exists a $T>0$ such that for every $\bar{\bm Q}$ satisfying equations~\eqref{eqn:fl_1}-\eqref{eqn:fl_3} with $\norm{\bar{\bm Q}(0)}_1\neq0$, we have
\begin{equation}
\bar{Q}_r(t)=0, \quad r\in \mathcal{R},\ \forall\, t\geq T.
\end{equation}
A fluid limit is said to be \emph{throughput optimal} if it is stable for all $\bm \lambda \in \mathcal C^\circ$.
\end{definition}

Thus far we have not proven that there is a direct link between an MDPN and its fluid model. However, if we treat the fluid model as a model in its own right, we can prove its throughput optimality as follows:
\begin{theorem}
\label{thm:fluid_stability}
The WARP fluid model (\ref{eqn:fl_1}-\ref{eqn:fl_3}) is throughput optimal.
\end{theorem}

\Beginproof{}
To prove the fluid stability, we show that there exists a $T>0$ such that for every fluid solution $\brac{\bar{\bm A},\bar{\bm Q},\bar{\bm D}}$ satisfying equations~\eqref{eqn:fl_1}-\eqref{eqn:fl_3} with $\norm{\bar{\bm Q}(0)}_1 = 1$, we have
$
\bar{Q}_r(t)=0, \ r\in \mathcal{R}, \forall \ t\geq T.
$
Consider a Lyapunov function $L(\bar{\bm Q}(t)) = \frac{1}{2}\sum_{r \in \mathcal R}\bar{Q}_r^2(t).$

Using fluid model equations~\eqref{eqn:fl_1}-\eqref{eqn:fl_3}, 
the derivative of $L(\bar{\bm Q}(t))$ is:
\begin{align}
    \frac{dL(\bar{\bm Q}(t))}{dt}&= \sum_{r\in \mathcal{R}}\bar{Q}_r(t) \lambda_r -\sum_{r\in \mathcal{R}} \bar{Q}_r(t) \bar{D}'_r(t)  
    \label{eqn:lyapunov_drift}.
\end{align}
Now from Theorem~\ref{thrm:stabilityRegion}, we know that for any $\bm \lambda \in \mathcal{C}^\circ$ there exists an agnostic policy $\pi'$ such that for some $\epsilon \in (0,1)$ we have 
\begin{equation}
\label{eqn:int_1}
\lambda_r + \epsilon  < 
\E_{\bm z \sim \mu^{{\pi'}}}[\sigma^{\pi'}_r(\bm z)]\, ,   .
\end{equation}
$\forall r  \in \mathcal R \,$.
Using~\eqref{eqn:int_1} in~\eqref{eqn:lyapunov_drift} we can write 
\begin{align}
    \frac{dL(\bar{\bm Q}(t))}{dt} 
    &\leq
    \sum_{r \in \mathcal R}\bar{ Q}_r(t)  \E_{\bm z \sim \mu^{{\pi'}}}[\sigma^{\pi'}_r(\bm z)]\notag
    \\
    &
    \quad
    -\sum_{r\in \mathcal{R}} \bar{Q}_r(t) \bar{D}'_r(t)
    -\epsilon
    \sum_{r \in \mathcal R}\bar{ Q}_r(t).
    \label{eqn:lyapunov_drift_2}
\end{align}
Note that the term $\sum_{r\in \mathcal{R}}\sum_{\bm z \in \mathcal Z}\mu^{\pi'}(\bm z)\sum_{a \in \mathcal{A}}p(a | \bm z)   \sigma_r \bar{Q}_r(t)$ corresponds to the average reward under the policy $\pi'$. 
Integrating and applying \eqref{eqn:fl_3} and then~\eqref{eqn:int_1} we see that 
\begin{align*}
&L(\bar{\bm Q}(t))
-
L(\bar{\bm Q}(s))
\\ 
&=\int_{s}^t 
    \sum_{r \in \mathcal R}\bar{ Q}_r(u)  \E_{\bm z \sim \mu^{{\pi'}}}[\sigma^{\pi'}_r(\bm z)]du
\\
&
\leq 
\int_{s}^t 
    \sum_{r \in \mathcal R}\bar{ Q}_r(u)  \E_{\bm z \sim \mu^{{\pi'}}}[\sigma^{\pi'}_r(\bm z)]
du
\\ & \qquad
    -
\int_{s}^t 
    \sum_{r\in \mathcal{R}} \bar{Q}_r(u) d\bar{D}_r(u)
    -
\int_{s}^t 
    \epsilon
    \sum_{r \in \mathcal R}\bar{ Q}_r(u)   
du
\\
&
\leq 
    -
\int_{s}^t 
    \epsilon
    \sum_{r \in \mathcal R}\bar{ Q}_r(u)   du
\\
&
\leq -\frac{\epsilon }{\sqrt{|\mathcal R|}}\int_s^t L(Q(u))^{\frac{1}{2}}du \, .
\end{align*}
The last inequality follows using the bound: $\norm{\bar{\bm Q}}_1\geq {\norm{\bar{\bm Q}}_2} = (2 L(\bar{\bm Q}))^{\frac{1}{2}}$. From this, we see that at any point of differentiability
\begin{equation}
\label{eqn:lyapunov_drift_3}
\frac{dL(\bar{\bm Q}(t))}{dt} \leq -{\epsilon } \brac{2 L(\bar{\bm Q}(t))}^{1/2}.   
\end{equation} 
From~\eqref{eqn:lyapunov_drift_3} it can be observed that if $L(\bar{\bm Q}(T)) = 0$ at any differentiable point $ T$, then $L(\bar{\bm Q}(t)) = 0 $ for all $t \geq T $. On the other hand, while $L(\bar{\bm Q}(t)) > 0$, we have from~\eqref{eqn:lyapunov_drift_3} 
\begin{align*}
&\brac{L(\bar{\bm Q}(t))}^{1/2}-\brac{L(\bar{\bm Q}(0))}^{1/2}
\\
&=\frac{1}{2}\int_0^t \brac{2 L(\bar{\bm Q}(t))}^{-1/2} \frac{dL(\bar{\bm Q}(t))}{dt} dt \leq -\frac{\epsilon }{\sqrt{2}}t.
\end{align*}
Thus,
$
L(\bar{\bm Q}(t)) \leq ( 2 L(\bar{\bm Q}(0))^{1/2} -{\epsilon t }/\sqrt{2} )_+^2.
$
Moreover, the function $L(\bar{\bm Q}(t))$ is continuous and non-increasing. Hence, for $\norm{\bar{\bm Q}(0)}_1\neq0$, $L(\bar{\bm Q}(t))=0$ for all $t\geq T$ where 
$
T=
{
        2\sqrt{|\mathcal R|}\brac{L(\bar{\bm Q}(0))}^{1/2}
    }/{
        \epsilon
        },
$
and we have $\bar{Q}_r(t)=0$ for $r\in \mathcal{R}$, for all $t\geq T$. This completes the proof of Theorem~\ref{thm:fluid_stability}. 
\Endproof

\subsection{Fluid Limit and Timescale Separation}
\label{sec:timescale}

In this section, we first introduce the fluid scaling of the system processes. We then present the fluid limit and discuss the concept of timescale separation, which is crucial for analyzing the system's asymptotic behavior in the fluid limit.

\subsubsection{Fluid Scaling}

We define the fluid scaled processes $\bar{\bm X}^c=\brac{\bar{\bm X}^c(t): t\in \Z_+}$ where $\bar{\bm X}^c(t)=\brac{\bar{\bm Q}^c(t),\bar{\bm Z}^c(t)}$  such that $\forall t\in \Z_+, \ c\in \Nats $
\begin{equation*}
\bar{\bm Q}^c(t)=\frac{\mf Q(ct)}{c}, \qquad  \bar{\bm Z}^c(t)=\frac{\mf Z(ct)}{c}. 
\end{equation*}

Here, we accelerate time by the scaling factor $c$ and scale space by this same factor.
An important facet is the two-time scale separation between the processes $\bar{\bm Z}^c$  and $\bar{\bm Q}^c$ as discussed in Section \ref{sec:timescale}. 
In the limit as $c \to \infty$, the $\bar{\bm Z}^c$  process is zero, yet still has a meaningful impact on the evolution of the process $\bar{\bm Q}^c$.
Further in the sequence above, we allow the arrival rate $\bm \lambda^c$ to depend on $c$.

\subsubsection{Fluid Limit}

The following result is somewhat technical. It establishes that limits exist for sequences $\brac{\bar{\bm Q}^c}_{c}$, and when they exist, it proves that they must satisfy the fluid equations \eqref{eqn:fl_1}-\eqref{eqn:fl_3}. Thus, this proves that our fluid limit model equations represent the asymptotic behavior of our queueing process. 

\begin{theorem}[Fluid Limit]
\label{thm:fluid_convergence}
Given Assumption \ref{ass:tau} and $\bar{\bm Q}(0) \in \R_+^{|\mathcal{R}|}$ with $ \norm{\bar{\bm Q}^c(0)}_1 \leq 1 $ and $ \bar{\bm Q}^c(0){\to} \bar{\bm Q}(0) $ a.s. and in $\ell_1$ and $\bm \lambda^c \rightarrow \bm \lambda$  as $c \rightarrow \infty$. Then the sequence of stochastic processes $\brac{\bar{\bm{Q}}^c}_{c \in \mathbb{N}}$ under the WARP policy is tight with respect to the topology of uniform convergence on compact time intervals. Additionally, every weakly convergent subsequence of $\brac{\bar{\bm{Q}}^c}_{c \in \mathbb{N}}$ converges to a Lipschitz continuous process $\bar{\bm Q}$. This limiting process satisfies the fluid model equations~\eqref{eqn:fl_1}-\eqref{eqn:fl_3}.
\end{theorem}
Hence, the process $\bar{\bm Q}$, defined in Theorem~\ref{thm:fluid_convergence}, characterizes the dynamics of an MDPN under the WARP policy in the limit as $c \to \infty$.

It is important to note that the scheduling schemes for MDPNs or classical switches can respond relatively quickly to changes in request queue lengths. By setting $\tau(c)$ such that Assumption \ref{ass:tau} holds, the system has sufficient time for the server process to reach its steady state between changes in request queue lengths. 
Setting $\tau$ to grow as somewhat modestly with the queue size ensures that request queue updates are not too frequently affecting server process and the queue size process maintains efficiency.

In our following result, we show that for the fluid limit solutions of every weakly convergent subsequence $\brac{\bar{\bm{Q}}^c}_{c \in \Nats}$, the schedule $\bm \sigma^\star$ is selected by the policy $\bar{\pi}^\star$, which is the limit of the sequence of policies $\brac{\pi^\star}_{c \in \Nats}$, and is obtained from~\eqref{eq:AREMDP}.

\subsubsection{Timescale Separation}

As the system becomes congested and the queue sizes grow large, the server dynamics operate on a significantly shorter timescale compared to the request process. This separation allows us to approximate the service process by its stationary behavior when analyzing the long-term evolution of the queues. This is a non-standard methodological step required to establish the fluid limit and connection between average reward MDPs and throughput optimality results for MDPNs.

\begin{restatable}{proposition}{TimeScale}\label{prop:timescale}
 For any interval $[u,t]$ for which $\bar Q_r(s)>0$ for all $s \in [u,t]$ the following holds
\begin{equation}
\bar{D}_r(t) - \bar{D}_r(u)  = \lim_{c\rightarrow \infty} \frac{1}{c} \sum_{s=cu}^{ct} \E_{z\sim \mu^{\pi^\star}(s)} \sbrac{\sigma^\star_r }\,,
\label{eq:Dlim}    
\end{equation}
where above $\mu^{\pi^\star}(s)$ is the stationary distribution of the servers under the WARP policy at time $s$ and $\sigma_r^\star$ is the stationary number of type $r$ requests under $\mu^{\pi^\star}(s)$.   
\end{restatable}

\subsection{Stochastic Stability.}

With Proposition \ref{prop:timescale} and Theorems \ref{thm:fluid_stability} and \ref{thm:fluid_convergence} in place, the proof of  Theorem \ref{thm:stochastic_stability}
is a technical albiet somewhat standard application of the Multiplicative Foster-Lyapunov Theorem (see \cite{robert2013stochastic} and \cite{bramson2008stability}). The details of this final step are given in Section \ref{sec:fluid_limit_proof} of the appendix.


\section{Conclusions}
This work proposes an optimal scheduling scheme for decision dependent servers with the Markov Decision Processing Network (MDPN) model. We analyze a general MDP formulation of service, which is significantly different from the simplified service networks and topologies explored in prior works.  
A key contribution is a novel, throughput-optimal policy called the WARP policy. The WARP policy optimizes server utilization to serve incoming requests efficiently. Furthermore, the paper makes significant technical contributions by developing a novel method to establish the fluid limit using two-time scale separation in a general MDPN. Our work has substantive implications for advances in queueing networks and related applications. A future direction is to establish guarantees for the multiobjective reinforcement learning algorithms in the context of these queueing networks.

\bibliographystyle{plainnat} 
\bibliography{OR_Switch}

\appendix

\section{Proof of Theorem \ref{thrm:stabilityRegion}}
\label{append:stabilityregion}

We now restate and prove Theorem \ref{thrm:stabilityRegion}.

\ThrmStabilityRegion*

\Beginproof{}
We first prove \eqref{eq:Ccirc} holds. To do this, we take any stabilizable policy. We then use the long-run service process of that policy to design a server scheduling scheme with the correct stationary service rate to satisfy \eqref{eq:Ccirc}.

First, if $\bm \lambda \in \mathcal C$, then let $\bm X(\infty)=(\bm Z(\infty),\bm Q(\infty))$ denote the stationary distribution of the Markov chain that is positive recurrent under $\bm \lambda$. Also, we let $ d_r(\infty)$ denote the stationary number type $r$ requests that depart. Thus, we have for $r \in \mathcal R$,
\begin{align}
    \lambda_r 
    &\leq \mathbb{E}\left[ d_r(\infty) \right] \notag \\
    &= \sum_{\bm z \in \mathcal Z} \sum_{\bm q \in \mathcal Q} 
        \mathbb{P}\left( \bm X(\infty) = (\bm z, \bm q) \right)
        \sum_{a \in \mathcal{A}} 
            P(a \mid \bm z, \bm q) \, {\sigma}_r(a, \bm z)
    \label{eq:expand}
\end{align}
Above, \eqref{eq:expand} the expanded expression for the stationary departure rate.
Here $P(a \mid \bm z, \bm q) $ denotes the stationary probability that action $a$ is selected given the system state $(\bm z, \bm q)$. 

Conditioning on $\bm Z(\infty)$, we can further write
\begin{align}
 \lambda_r
&\leq\sum_{\bm z \in \mathcal Z}
    \mb{P}(\bm Z(\infty)=\bm z) 
    \sum_{\bm q \in \mathcal Q}
        \mb{P}(\bm Q(\infty)=\bm q | \bm Z(\infty) =  \bm z)
    \sum_{a\in\mathcal{A}}
            p(a | \bm z, \bm q)
             {\sigma}_r(a,\bm z)
            \notag \\
 & =\sum_{\bm z \in \mathcal Z}\mb{P}(\bm Z(\infty)=\bm z)\sum_{a \in \mathcal{A}}p(a | \bm z)   {\sigma}_r(a,\bm z) \label{eq:lle_stationary}
\end{align}
where 
\begin{equation}\label{eq:pidef}
p(a | \bm z)
:=
\sum_{\bm q \in \mathcal Q}
    \mb{P}(\bm Q(\infty)=\bm q\mid \bm Z(\infty)=\bm z)
    p(a | \bm z, \bm q) 
\end{equation}
 denotes the probability that schedule $a$ is selected in stationary regime given that $\bm Z(\infty)=\bm{z}$. 

We can now use it to define an entanglement matching policy that is agnostic to request queues. In particular, we define the request agnostic policy $\hat \pi$ where whenever the left side queues are in state $\bm z$, we choose the action $a$ with probability $p(a | \bm z)$ as given by \eqref{eq:pidef}. We let $\hat{\bm Z}(\infty)$ denote the stationary number of servers under this policy. 

We now show that $\hat{\bm Z}(\infty)$ and $\bm Z(\infty)$ are equal in distribution.\footnote{For a standard switching model these are immediately equal as they are both independent of the queue size process $\bm Q(\infty)$; however, in our case this property that must be verified. The following calculations are not standard in prior analysis on Maximal Stability.} We prove this by verifying that both distributions satisfy identical balance equations.

First the balance equations for $\hat{\bm Z}(\infty)$  are
\begin{equation}
\mathbb P (\hat{\bm Z}(\infty) = \bm z) 
=
\sum_{\bm z' \in \mathcal Z} \sum_{a \in \mathcal{A}} \mathbb P 
( \bm Z(1) = \bm z | \bm Z(0) =\bm z', a(0) = a ) p(a | \bm z')
\mathbb P( \hat {\bm Z}(\infty) = \bm z' ) \, .
\label{eq:hatbalance}    
\end{equation}
Here, $\bm z'$ indicates the initial state of the stationary distribution, and $\bm z$ represents the next state reached after one transition.
Next, the balance equations for $\bm X(\infty) = (\bm Z(\infty), \bm Q(\infty))$ are 
:
\begin{align}
    &
    \mathbb P (  \bm Z( \infty) = \bm z, \bm Q(\infty) = \bm q )
    \notag
\\
&=
    \sum_{\bm z' \in \mathcal Z} \sum_{\bm q' \in \mathcal Q} \sum_{a \in \mathcal{A}} 
    \Big[
    \mathbb P ( \bm Z(1) =\bm z  , \bm Q(1) = \bm q  | \bm Z(0) =\bm z ' , \bm Q(0) = \bm q', a(0)= a)
    \notag
    \\
    &\qquad \qquad \qquad \qquad
    \times 
    p(a | \bm z', \bm q')
    \notag
    \\
    &\qquad \qquad \qquad  \qquad\qquad
    \times
    \mathbb P (\bm Z(\infty) =\bm z ' , \bm Q(\infty ) = \bm q' )  
    \Big]\, .
     \label{eq:Zlong0}
\end{align}
Thus, summing over $\bm q$ gives
\begin{align}
  &\mathbb P (  \bm Z( \infty) = \bm z)
  \notag
  \\
&
  = \sum_{\bm q \in\mathcal Q} \mathbb P (  \bm Z( \infty) = \bm z, \bm Q(\infty) = \bm q ) 
  \notag
  \\
  &
  =
    \sum_{\bm q, \bm z' \bm q'} \sum_{a \in \mathcal{A}} 
    \Big[ 
    \mathbb P ( \bm Z(1) =\bm z  , \bm Q(1) = \bm q  | \bm Z(0) =\bm z ' , \bm Q(0) = \bm q' , a(0)= a)  
    \notag
    \\
    &
    \qquad \qquad \qquad \qquad 
    \times
    p(a | \bm z', \bm q')
    \notag
    \\
    &\qquad \qquad \qquad \qquad\quad \times
    \mathbb P (\bm Z(\infty) =\bm z ' , \bm Q(\infty ) = \bm q' )  \Big]
    \label{eq:Zlong1}
\\
&
=
    \sum_{\, \bm z' \in \mathcal Z, \bm q' \in \mathcal Q} \sum_{a \in \mathcal{A}} 
   \Big[ 
    \mathbb P ( \bm Z(1) =\bm z  | \bm Z(0) =\bm z ' , \bm Q(0) = \bm q' , a(0)= a )  
    \notag
    \\
    &
    \qquad \qquad \qquad \qquad 
    \times
    p( a | \bm z', \bm q')
    \notag
    \\
    &\qquad \qquad \qquad \qquad\quad \times
    \mathbb P (\bm Z(\infty) =\bm z ' , \bm Q(\infty ) = \bm q' )  \Big]
    \label{eq:Zlong2}
\\
&
=
    \sum_{\, \bm z' \in \mathcal Z} \sum_{\bm q' \in \mathcal Q} \sum_{a \in \mathcal{A}} 
   \Big[ 
    \mathbb P ( \bm Z(1) =\bm z  | \bm Z(0) =\bm z ' ,  a(0)= a) 
    \notag
    \\
    &
    \qquad \qquad \qquad \qquad 
    \times
    p(a| \bm z', \bm q')
    \notag
    \\
    &\qquad \qquad \qquad \qquad\quad \times
    \mathbb P ( \bm Q(\infty ) = \bm q' 
 | \bm Z(\infty) =\bm z '  ) 
     \mathbb P (\bm Z(\infty) =\bm z ' ) \Big]  
     \label{eq:Zlong3}
\\
&
=
    \sum_{\, \bm z' \in \mathcal Z}
    \sum_{a \in \mathcal{A}} 
    \mathbb P ( \bm Z(1) =\bm z  | \bm Z(0) =\bm z ' ,  a(0)= a ) 
    \notag
    \\
    &
    \qquad \qquad \qquad  
    \times
    \Big[  
    \sum_{\, \bm q' \in \mathcal Q}
   p(a | \bm z', \bm q')
    \mathbb P ( \bm Q(\infty ) = \bm q' 
 | \bm Z(\infty) =\bm z '  ) 
 \Big]
     \mathbb P (\bm Z(\infty) =\bm z ' )   
      \label{eq:Zlong4}
\\
&
=
    \sum_{\, \bm z' \in \mathcal Z}
    \sum_{a \in \mathcal{A}} 
    \mathbb P ( \bm Z(1) =\bm z  | \bm Z(0) =\bm z ' ,  a(0)= a )
    p(a | \bm z')
         \mathbb P (\bm Z(\infty) =\bm z ' )   \,.
          \label{eq:Zlong5}
\end{align}
Above in \eqref{eq:Zlong1}, we apply \eqref{eq:Zlong0}. In \eqref{eq:Zlong2} we sum over $\bm q$. In \eqref{eq:Zlong3}, we note that the transition from $\bm z'$ to $\bm z$ does not depend on $\bm q'$ once we condition on the action $\bm a$, and we also condition on $\bm z'$. In Equation \eqref{eq:Zlong4}, we notice we can bring the sum over $\bm q'$ inside. In Equation \eqref{eq:Zlong5}, we note by our definition \eqref{eq:pidef} that this inner sum is $p(a | \bm z')$. From \eqref{eq:Zlong5} and \eqref{eq:hatbalance} we see that both $\bm Z(\infty)$, $\hat {\bm Z}(\infty)$ satisfy the same balance equation (on the same irreducible set of states). Thus, by the uniqueness of stationary distributions $\bm Z(\infty)$, $\hat {\bm Z}(\infty)$ are equal in distribution.

We can now see that \eqref{eq:Ccirc2} holds. In particular, we let $\mu (\bm z) = \mathbb P(\hat {\bm Z}(\infty) = \bm z) = \mathbb P(\bm Z(\infty) = \bm z)$ then, from \eqref{eq:expand}, we see that our request agnostic policy $\hat \pi$ is such that 
\[
\lambda_r
 \leq \sum_{\bm z \in \mathcal Z }\mu(\bm z)\sum_{a \in \mathcal{A}}p(a | \bm z)   {\sigma}_r(a,\bm z)   \, r \in \mathcal R
\]
This verifies \eqref{eq:Ccirc2} and completes the first part of the theorem.

For the second part, we apply Foster's Lemma.\footnote{Again, there are some complications when compared to more traditional maximal stability proofs due to the random evolution of the service process.} Consider an agnostic policy such that 
\[
\lambda_r
 < \sum_{\bm z \in \mathcal Z }\mu(\bm z)\sum_{a \in \mathcal{A}}p(a | \bm z)   {\sigma}_r(a,\bm z)  \,, \quad r \in \mathcal R \, ,
\]
holds.
Notice that if the server process is stationary, then for all $\bm q$ sufficiently large,
\begin{align*}
    &\mathbb{E}\left[
        \sum_{r \in \mathcal{R}} Q_r(t+1)
        - \sum_{r \in \mathcal{R}} Q_r(t)
        \,\Big|\, \bm Q(t) = \bm q,\; \bm Z(t) \sim \mu
    \right]
    \\
    &= \sum_{r \in \mathcal{R}}
        \left(
            \lambda_r
            - \sum_{\bm z \in \mathcal{Z}} \mu(\bm z)
                \sum_{a \in \mathcal{A}} p(a \mid \bm z)\, {\sigma}_r(a, \bm z)
        \right)
    \\
    &< 0 \, .
\end{align*}
Thus, in principle, if it were not for Markov evolution of the $\bm Z(t)$ process, then we could directly apply Foster's Lemma. To address this issue, we simply need to allow the process sufficient time to approach equilibrium, and then we can apply Lyapunov arguments. We begin by proving stability for a single queue, as $(\bm Z(t), Q_r(t))$ is a Markov chain for agnostic policies. We then use our queue size bound to prove stability for $(\bm Z(t),  \bm Q(t))$. Furthermore, we apply the Lyapunov argument of \cite{bertsimas1998geometric}. This is the plan for the remainder of the proof. \footnote{We note that an alternative approach is to apply Foster-Lyapunov over renewal cycles of the process $\bm Z$ and apply Theorem 8.13 of \cite{robert2013stochastic}. Another approach would be to consider a quadratic Lyapunov function and then apply a similar argument to \cite{georgiadis2006resource}.}

We let 
\[
\epsilon_r :=  \Big( \sum_{\bm z}\mu(\bm z)\sum_{a \in \mathcal{A}}p(a | \bm z)   {\sigma}_r(a,\bm z) \Big)   -\lambda_r  \, .
\]
We note that under this policy, the lefthand process $\bm Z(t)$ is an ergodic Markov chain thus there exists a $t_0$ such for all $t\geq t_0$ 
\[
\left|\sum_{\bm z}\mu(\bm z)\sum_{a \in \mathcal{A}}p(a | \bm z)   {\sigma}_r(a,\bm z)  -
\sum_{\bm z}\mathbb P (\bm Z(t) = \bm z)\sum_{a \in \mathcal{A}}p(a | \bm z)   {\sigma}_r(a,\bm z)  \right| < \frac{\epsilon_r}{2}\, .
\] 

\Endproof

\section{Stochastic Stability: Proof of Theorem \ref{thm:stochastic_stability}}

\subsection{Fluid Limit: Proof of Theorem~\ref{thm:fluid_convergence}}
\label{sec:fluid_limit_proof}
The proof of Theorem~\ref{thm:fluid_convergence} is structured as follows.
We show that the sequence of stochastic processes $\brac{\bar{\bm{Q}}^c}_{c \in \Nats}$ has a subsequence with a limit. This argument is standard -- albeit somewhat technical -- and employs the Arzelà-Ascoli (See Theorem~\cite{billingsley2013convergence}).
With the existence of limits confirmed, we explore their properties.
In particular, in Proposition \ref{prop:timescale},
 we examine the timescale separation of $\bar{\bm{Q}}^c$ and $\bar{\bm{Z}}^c$ over the time interval $ct$ by partitioning time into time intervals of length $\tau(\bar{\bm{Q}}^c)$. 
  Within each length $\tau(\bar{\bm{Q}}^c)$ interval, the actions are selected by policy $\pi^\star$ are determined based on the request queue lengths at the start of the segment. We can then apply this in Proposition \ref{prop:opt} to give fluid equation \eqref{eqn:fl_3}.

\subsubsection{Tightness.}
\label{subsec:tightness}
We first define the family of coupled processes $\brac{\bar{\bm{Q}}^c, \bar{\bm{A}}^c, \bar{\bm{D}}^c }$, where $\bar{\bm{Q}}^c(t)$, $\bar{\bm{A}}^c(t)$, and $\bar{\bm{D}}^c(t)$ are the scaled versions of the request queue process, the cumulative arrivals and departures, respectively. These processes are constructed on the same probability space and remain the same for different values of $c$.
For each $r \in \mathcal{R}$ and $l\in \mathcal{Z}$, we define the scaled processes for arrivals, departures, and server queues as:
\begin{equation}
\bar{A}_r^c(t)=\frac{ A_r(ct)}{c}, \;\bar{D}_r^c(t)=\frac{ D_r(ct)}{c},\; \bar{Z}_l^c(t)=\frac{Z_l(ct)}{c} .
\end{equation}
for $c\in \Nats$.

The proof of tightness is contained in the following proposition.

\begin{restatable}{proposition}{PropositionTightness}\label{prop:Tightness}
       The sequence of stochastic processes $\brac{\bm \lambda^c,\bar{\bm{Q}}^c,\bar{\bm{A}}^c,\bar{\bm{D}}^c, \bar{\bm Z}^c}_{c \in \mathbb{N}}$ under the WARP policy is tight with respect to the topology of uniform convergence on compact time intervals.  
\end{restatable}

\Beginproof{}
To prove tightness, we show that there exists a measurable set $G$ with $\mathbb{P}(G) = 1$, such that for all $\omega \in G$, any subsequence of $\brac{\bm \lambda^c,\bar{\bm{Q}}^c,\bar{\bm{A}}^c,\bar{\bm{D}}^c, \bar{\bm Z}^c}_{c \in \mathbb{N}}$ under the WARP policy contains a further subsequence that converges uniformly on compact time intervals.

Firstly, since $\bm \lambda^c$ belongs to the compact set $\mathcal C$, we chose a subsequent along which $\bm \lambda^c$ converges to some value $\bm \lambda$. 
As discussed in Subsection~\ref{subsec:topology}, for each $r \in \mathcal{R}$, the collection $(A_r(t) - A_r(t-1) : t \in \Nats)$ consists of i.i.d. random variables with mean $\lambda_r$. Therefore, by the Functional Strong Law of Large Numbers, on a set $G_1$ with $\mathbb{P}(G_1) = 1$, we have for each $r \in \mathcal{R}$ 
$$\bar{A}^c_r(t) \to \lambda_r t,$$
as $c \to \infty$, with the convergence being uniform on compact intervals. Similarly, since $Z_l(t)$ is bounded we have
\[
\bar{Z}_l^c(t) \rightarrow 0
\]
as $c\rightarrow\infty$, with the convergence being uniform on compact intervals.
From the Arzelà-Ascoli Theorem, we know that any sequence of equicontinuous functions $\bar{Y}^c(t)$ on $[0,T]$ for $T > 0$, with $\sup_c\left|\bar{Y}^c(0)\right| < \infty$, has a converging subsequence with respect to the uniform norm. We will now verify that any subsequence of $\brac{\bm \lambda^c,\bar{\bm{Q}}^c,\bar{\bm{A}}^c,\bar{\bm{D}}^c, \bar{\bm Z}^c}_{c \in \mathbb{N}}$ satisfies both conditions for any $\omega \in G_1$.

Since $\bar{A}_r^c(0)$ and $\bar{D}_r^c(0)$ are initially $0$ and $\|\bar{\bm Q}^c(0)\|_1 \leq 1$, the supremum is also bounded. The equicontinuity of the sequence $\bar{A}_r^c$ follows from uniform convergence on compact intervals. Moreover, the fact that the maximum number of requests served in each time slot is bounded above by $B$ implies that $\bar{D}_r^c$ is Lipschitz continuous, which in turn implies the equicontinuity of this sequence. Finally, since $\bar{Q}_r^c$ is the sum of a bounded number of equicontinuous functions, the sequence $\bar{Q}_r^c$ is also equicontinuous. Therefore, as the conditions of the Arzelà-Ascoli Theorem are satisfied for all $\omega \in G_1$, every subsequence of $\brac{\bm \lambda^c,\bar{\bm{Q}}^c,\bar{\bm{A}}^c,\bar{\bm{D}}^c, \bar{\bm Z}^c}_{c \in \mathbb{N}}$ contains a further subsequence that converges uniformly. Moreover, since these sequences of functions are uniformly Lipschitz continuous, their limits must also be Lipschitz continuous.
\Endproof

Tightness implies relative compactness \cite[Prohorov's Theorem]{billingsley2013convergence}: that for every $\brac{\bm \lambda^c,\bar{\bm{Q}}^c,\bar{\bm{A}}^c,\bar{\bm{D}}^c, \bar{\bm Z}^c}_{c \in \mathbb{N}}$ there is a weakly convergent sub-sequence.
This verifies the tightness statement in Theorem \ref{thm:fluid_convergence} and also the fluid limit equation \eqref{eqn:fl_1}.
Applying the Skorohod Representation Theorem \cite{billingsley2013convergence} we may also assume that this sub-sequence convergence holds almost surely. 
Throughout the rest of the proof we assume that $\brac{\bm \lambda^c,\bar{\bm{Q}}^c,\bar{\bm{A}}^c,\bar{\bm{D}}^c, \bar{\bm Z}^c}_c$ is a sub-sequence that converges almost surely uniformly on compacts to its limit point is $\brac{\bm \lambda,\bar{\bm{Q}},\bar{\bm{A}},\bar{\bm{D}}, \bar{\bm Z}}\,$. By applying the Functional Strong Law of Large Numbers, $\bar{A}_r^c(t)$ converges uniformly on compacts to $ \lambda_r t$  as $c \to \infty$ this verifies \eqref{eqn:fl_2}. It remains to prove \eqref{eqn:fl_3}.

\subsubsection{Timescale Separation.}
To prove \eqref{eqn:fl_3}, we first need to understand the timescale separation behavior of these switches on the fluid scale. The main finding of this section is summarized in the following proposition.

\TimeScale*

\Beginproof{}
To simplify notation in the proof, without loss of generality, we assume that $u=0$.
Since $\bar{Q}_r(s)>0$ for $s\in [0,t]$ and $\bar{Q}^c_r(s)$ convergences uniformly to $\bar{Q}_r(s)>0$, there exists a value of $c'$ and $\epsilon >0$, for which $\bar{Q}^c_r(s) \geq \epsilon$ for all $s \in [0,t]$,  for all $c\geq c'$. Thus, since all queues are non-zero all scheduled requests are served, and so
we can write $\bar{D}^c_r(t)$ as
\begin{align}
\bar{D}^c_r(t)&= \frac{1}{c}\sum_{s=1}^{ct}\sigma^\star_r\brac{\bm Z(s);\bar{\bm Q}^c(s/c)}\label{eqn:scales_Q}\, .
\end{align}
Also note that the policy $\pi^\star(\cdot)$ and resulting schedule $\sigma^\star(\cdot)$ are obtained as a solution of the average reward MDP, and these solutions do not change if we rescale the queue size vector by a factor $c$; that is, we have
$\pi^\star\brac{\cdot\; ;{\bm Q}(s)} \equiv \pi^\star({\cdot\;;\bar{\bm Q}^c(s/c)})$
and $\sigma^\star \brac{\cdot\; ;{\bm Q}(s)} \equiv \sigma^\star({\cdot\;;\bar{\bm Q}^c(s/c)})$.

Because the WARP policy fixes the policy over variable time intervals, we also introduce some notation to simplify the exposition. Henceforth, we denote the  beginning of each time interval with
\[
t(i+1) := t(i) + \tau\big( \bm Q(t(i))  \big) \, ,
\]
for epochs $i \in \mathbb N$ and with $t(0):=0$.
We also let $\pi_i^\star$ be the policy used by the WARP policy over the $i$th epoch. Specifically,
\[
\pi_i^\star (\bm z) := \pi^\star(\bm z;\bm Q(t(i))), 
\]
{for  states $\bm z \in \mathcal{Z}$}. Also, we define
\[
\sigma^\star_{r,i}(\bm z) := \sigma_r\big( \pi_i^\star(\bm z),\bm z;\bar{\bm Q}^c(t(i)/c)\big) \, .
\]
{for  states $\bm z \in \mathcal{Z}$}.

Now decompose the departure term, equation~\eqref{eqn:scales_Q}, over the intervals $[t(i), t(i+1))$. This decomposition is beneficial as it leverages the timescale separation between the server and request processes. More precisely, the departure term can be written as :
\begin{align}
    \bar{D}^c_r(t)
    &=
    \frac{1}{c} \sum_{i:t(i) \leq ct}
    \Bigg[
        \sum_{s=t(i)}^{t(i+1)-1}
        \sigma^\star_{r,i}\big(\bm Z(s)\big)
        -
        \sum_{s=t(i)}^{t(i+1)-1}
        \mathbb{E}\Big[
            \sigma^\star_{r,i}\big(\bm Z(s)\big)
            \,\Big|\, \bm X(t(i))
        \Big]
    \Bigg]
    \label{eqn:term1}
    \\
    &\quad
    + \frac{1}{c} \sum_{i:t(i) \leq ct}
    \Bigg[
        \sum_{s=t(i)}^{t(i+1)-1}
        \Big(
            \mathbb{E}\Big[
                \sigma^\star_{r,i}\big(\bm Z(s)\big)
                \,\Big|\, \bm X(t(i))
            \Big] \notag
            \\ & \qquad \qquad
            -
            \mathbb{E}\Big[
                \sigma^\star_{r,i}\big(\bm Z(s)\big)
                \,\Big|\, \bm Z(t(i)) \sim \mu^{\pi^\star_i},\, \bm Q(t(i))
            \Big]
        \Big)
    \Bigg]
    \label{eqn:term2}
        \\
    &\quad
    +  \frac{1}{c}\sum_{s=t(i^*)+1}^{t(i^*+1)-1}
        \sigma^\star_{r,i^*}\big(\bm Z(s)\big)
        \label{eqn:term4}
    \\
    &\quad
    + \frac{1}{c} \sum_{i:t(i) \leq ct}
    \sum_{s=t(i)}^{t(i+1)-1}
    \mathbb{E}\Big[
        \sigma^\star_{r,i}\big(\bm Z(\infty)\big)
        \,\Big|\, \bm Z(t(i)) \sim \mu^{\pi^\star_i},\, \bm Q(t(i))
    \Big]
    \label{eqn:term3}
\end{align}
where above $i^* = \max\{ i : t(i) \leq ct\}$. Above we note that the terms in \eqref{eqn:term2} and \eqref{eqn:term3} cancel; however, we emphasize the stationarity of the process $\bm Z$ in \eqref{eqn:term3}. 

The plan for the remainder of the proof is as follows: for term~\eqref{eqn:term1}, the Azuma-Hoeffding inequality will be used to show that it approaches zero as $c \to \infty$. 
Moreover, term~\eqref{eqn:term2} represents the difference between the expected departures when the process $\bm{Z}$ starts from its stationary distribution $\mu^{\pi^\star_i}$ under the policy $\pi^\star_i$ as defined by \eqref{eq:AREMDP} and when $\bm{Z}$ begins its actual state at time $t(i)$. Leveraging timescale separation will show that term~\eqref{eqn:term2} vanishes as $c \to \infty$.
The term \eqref{eqn:term4} is an error term due to the discretization over intervals $[t(i),t(i+1))$, we will see that this term vanishes since the length of intervals is always of order $o(c)$.
From this we will see that the convergence of $\bar{D}_r^c(t)$ term~\eqref{eqn:term3} is determined by the term ~\eqref{eqn:term3} as $c \to \infty$. So, we now analyze terms~\eqref{eqn:term1}-\eqref{eqn:term3} separately.

 First for term \eqref{eqn:term1}, we notice that 
\[
\sum_{s=t(i)}^{t(i+1)-1}
    \Bigg[
        \sigma^\star_{r,i}\big(\bm Z(t(i)+s)\big)
        -
        \mathbb{E}\big[\,\sigma^\star_{r,i}\big(\bm Z(t(i)+s)\big)
            \mid \bm X(t(i))\,\big]
    \Bigg] 
\]
is the sum of a martingale difference sequence. Therefore, by Lemma~\ref{lem:azzuma_hoeffding}, which applies the Azuma-Hoeffding inequality, the term of~\eqref{eqn:term1} converges to zero as $c\to \infty$.

Now, we turn our attention to term~\eqref{eqn:term2}. Here we first study the mixing time properties of the Markov chain $\bm Z$ under the fixed policy $\pi_i^\star$ over the time interval $[t(i), t(i+1))$. 

By \cite{rosenthal1995convergence}, we know that for any fixed policy $\pi$, since $Z$ is an irreducible and aperiodic Markov chain, there exist positive constants $C_\pi $ and $\rho_\pi <1$ such that for any initial distribution $\mu$ of $\bm Z(0)$, we have
\begin{equation}
\label{eq:mixing_bound}
\sup_{\mu}\| \mu [P^{\bm \pi}]^t - \mu^{\pi} \|_1\leq C_\pi t^{|\mathcal Z|} \rho_\pi^t,
\end{equation}
for $t \in \mathbb Z_+$. [The result follows by decomposing the transition matrix $P^{\bm \pi}$ into its Jordan normal form and expanding. The Perron-Frobenius theorem gives that $\rho_\pi<1$.] Since the number of policies is finite, we can choose $C' = \max_\pi C_\pi$ and $\rho = \max_\pi \rho_\pi <1$.
We now define the mixing time
\[
\tau_{\mix}(c) = \sbrac{\frac{|\mathcal Z|+2}{\log \rho}} \log c
\]
From this and the properties of the mixing time, it follows that:
\begin{align}
\sum_{\tau \geq \tau_{\mix}(c)} 
\sup_\pi 
\| 
\mu [P^{\bm \pi}]^t - \mu^{\pi} 
\|_1
    &\leq 
\sum_{\tau \geq \tau_{\mix}(c)} C \tau^{|\mathcal Z|} 
\rho^{\tau}
\notag
\\
    & \leq 
    C \tau_{\mix}(c)^{|\mathcal Z|}  \rho^{\tau_{\mix}(c)}
    \sum_{\tau \geq 0}
    \tau^{|\mathcal Z|} \rho^{\tau}
\notag
  \\
  &
  \leq
  C' \tau_{\mix}(c)^{|\mathcal Z|}  \rho^{\tau_{\mix}(c)} \frac{|\mathcal Z|!}{(1-\rho)^{|\mathcal Z|+1}}
  \notag
  \\
& \leq
        C' \brac{\frac{|\mathcal Z|+2}{\log \rho}}^{|\mathcal Z|}\frac{|\mathcal Z|!}{(1-\rho)^{|\mathcal Z|+1}} \frac{1}{c^2}
= \frac{C}{c^2} \, .
\label{eqn:mixing_bound_sum}
\end{align}
Above we define 
\[
C = C' 
    \left( 
        \frac{|\mathcal Z| + 2}{\log \rho} 
    \right)^{|\mathcal Z|} 
    \frac{|\mathcal Z|!}
         {(1 - \rho)^{|\mathcal Z| + 1}}
         \, .
\]




We can further decompose~\eqref{eqn:term2} as follows:
\begin{align}
\left|\eqref{eqn:term2}\right|
&\leq
\frac{1}{c}\sum_{i:\,t(i)\le ct}
    \left|
        \sum_{s=t(i)}^{t(i)+\tau_{\mix}(c)-1}
            \Big(
                \mathbb{E}\big[\sigma^\star_{r,i}(\bm Z(s))\mid \bm X(t(i))\big]
                -
                \mathbb{E}\big[\sigma^\star_{r,i}(\bm Z(s))\mid \bm Z(t(i))\sim \mu^{\pi_i^\star },\ \bm Q(t(i))\big]
            \Big)
    \right|
    \label{eqn:inter_1}
\\
&+
\frac{1}{c}\sum_{i:\,t(i)\le ct}
    \left|
        \sum_{s=t(i)+\tau_{\mix}(c)}^{t(i+1)-1}
            \Big(
                \mathbb{E}\big[\sigma^\star_{r,i}(\bm Z(s))\mid \bm X(t(i))\big]
                -
                \mathbb{E}\big[\sigma^\star_{r,i}(\bm Z(s))\mid \bm Z(t(i))\sim \mu^{\pi_i^\star },\ \bm Q(t(i))\big]
            \Big)
    \right|
    \label{eqn:inter_2}\,.
\end{align}
It is important to choose $\tau_{\mix}(c)$ so that the server process has (approximately) reached its steady state when the system reaches $\tau_{\mix}(c)$. Therefore, until $\tau_{\mix}(c)$, there is a need to sacrifice system performance. It will be shown that the payoff up to $\tau_{\mix}(c)$, corresponding to term~\eqref{eqn:inter_1}, is bounded. 
Using the triangle inequality and observing that the maximum number of requests served in a single time slot is bounded by 
$B$,~\eqref{eqn:inter_1} can be upper-bounded as follows:

\begin{align}
    &\frac{1}{c} \sum_{i: t(i) \le ct}
        \left|
            \sum_{s = t(i)}^{t(i) + \tau_{\mix}(c) - 1}
                \Big(
                    \E\sbrac{\sigma^\star_{r,i}(\bm Z(s)) \mid \bm X(t(i))}
                    -
                    \E\sbrac{\sigma^\star_{r,i}(\bm Z(s)) \mid \bm Z(t(i)) \sim \mu^{\pi^\star_i},\, \bm Q(t(i))}
                \Big)
        \right|
    \notag \\
    &\qquad\leq \frac{2 B \, \tau_{\mix}(c)\, t}{\min_{i: t(i) \le ct} \tau(\bm Q(t(i)))}
    \notag \\
    &\qquad\leq \frac{2 B \, \tau_{\mix}(c)\, t}{\tau(\epsilon c)}
    \xrightarrow[c\to \infty]{} 0
    \label{eqn:t_1_term_bound}
\end{align}
Above, we note that schedules are bounded, and we lower-bound the number of terms $i$ for which $t(i) \le ct$ holds [here we observe that $t(i+1)-t(i) \geq \tau(\bm Q(t(i))) \geq \min_{i: t(i) \le ct} \tau(\bm Q(t(i)))$] . We then note that $\bm Q(t(i)) \geq \bm 1 \epsilon c$ (as demonstrated in the first paragraph of this proof). The limit then follows by Assumption  \ref{ass:tau}.

Now for \eqref{eqn:inter_2}, we use the timescale separation to show that this term converges to zero as $c\to\infty$. By leveraging the assumption that policy $\pi^\star$ (and $\sigma^\star$) is fixed and thus $\bm Z$ evolves as a Markov chain during each time segment $\cbrac{t(i),\dots,t(i+1)}$. Using this, we can rewrite and upper bound~\eqref{eqn:inter_2} as follows:
\begin{align}
\eqref{eqn:inter_2}
&\le \frac{1}{c}
    \sum_{i:\,t(i) \le ct}
    \sum_{s=t(i)+\tau_{\mix}(c)}^{t(i+1)-1}
    \Bigg|
      \mu'_{\bm Q(t(i))}\big[P^{\bm Q(t(i))}\big]^s
         \bm \sigma^\star_{r,i}
      -
      \mu^{\pi^\star_i}\big[P^{\bm Q(t(i))}\big]^s
         \bm \sigma^\star_{r,i} 
    \Bigg|
\nonumber\\
&\le \frac{1}{c}
    \sum_{i:\,t(i) \le ct}
    \sum_{s=t(i)+\tau_{\mix}(c)}^{t(i+1)-1}
    2B \,
    \Big\|
      \mu'_{\bm Q(t(i))}\big[P^{\bm Q(t(i))}\big]^s
      -
      \mu^{\pi^\star_i}
    \Big\|_{\TV}
\label{eqn:inter_3}\\
&\le \frac{1}{c}
    \sum_{i:\,t(i) \le ct}
    \frac{2B C}{c^2}
\label{eqn:inter_4}\\
&\le \frac{2B C}{c^2} \xrightarrow[c\to \infty]{} 0
\label{eqn:inter_5}.
\end{align}

In the first inequality, $ \mu'_{\bm{Q}(t(i))}$ represents the initial distribution of the $\bm{Z}$ process, which may not be stationary.
Furthermore, \eqref{eqn:inter_3} is derived from: 
\begin{multline*}
 \left| \mu'_{\bm Q} \sbrac{P^{\bm Q}}^s \bm \sigma^\star_{r,i}-\mu^{\pi^\star_i}\sbrac{P^{\bm Q}}^s \bm \sigma^\star_{r,i} \right| 
 \leq B     \sum_{\bm z\in \mathcal{Z}} \left|\sum_{\bm z'\in \mathcal{Z}}\brac{\mu'_{\bm Q}(\bm z') -\mu^{\pi^\star_i}}\sbrac{P^{\bm Q}}^s_{\bm z' \bm z}\right|
 \\= 2B     \frac{1}{2}\sum_{\bm z\in \mathcal{Z}} \left|\mu'_{\bm Q}(\bm z)\sbrac{P^{\bm Q}}^s_{ \bm z} -\mu^{\pi^\star_i}(\bm z)\sbrac{P^{\bm Q}}^s_{ \bm z}\right|=2B \norm{\mu'_{\bm Q}\sbrac{P^{\bm Q}}^s-\mu^{\pi^\star_i}\sbrac{P^{\bm Q}}^s}_{\TV},
\end{multline*}
(For brevity, we suppress the dependence of $\bm{Q}$ on $ t(i)$). In~\eqref{eqn:inter_4}, we have used our mixing time bound \eqref{eqn:mixing_bound_sum}. 

Therefore, applying~\eqref{eqn:t_1_term_bound} and~\eqref{eqn:inter_5} to \eqref{eqn:term2}, we have 
\begin{equation*}
\left|\eqref{eqn:term2}\right|\leq \frac{2 B \tau_{\mix}(c)t}{\tau(\bm Q(t(i)))}+ \frac{4B\rho^{\tau_{\mix}(c)}t}{(1-\rho)\tau(\bm Q(t(i)))} \xrightarrow[c\to \infty]{} 0\, .
\end{equation*}

For the term \eqref{eqn:term4}, we have
\begin{align*}
\| \eqref{eqn:term4} \| \leq \frac{1}{c} B ( t(i^* +1) - t(i^*) ) \leq \frac{B}{c} \max_{0\leq s \leq c t} \tau(\bm Q(s))
\leq  \, \frac{B}{c} \tau(c + \lambda_{\max} ct) \xrightarrow[c \rightarrow \infty]{} 0 \, .
\end{align*}
In the first inequality, we use the fact that schedules are bounded by $B$. Then, in the second, the upper bound is by the largest epoch. We then note that queue sizes are bounded by the initial queue size (which, by definition, is $c$) and the number of arrivals, which is bounded above by $\lambda_{\max} ct$. Finally, by Assumption \ref{ass:tau}, we have that $\tau(c)/c \rightarrow 0$. 

Finally, we turn to the limit of~\eqref{eqn:term3}. Since we are staring $Z(t(i))\sim \mu^{\pi^\star(\bm Q(t(i))}$, we can write
\begin{align*}
&\frac{1}{c} \sum_{i : t(i) \leq ct} \sum_{s = t(i)}^{t(i+1)-1}
    \mathbb{E}\left[
                \sigma^\star_{r,i}\big(\bm Z(\infty)\big)
        \,\Big|\, \bm Z(t(i)) \sim \mu^{\pi^\star_i},\, \bm Q(t(i))
    \right]
\\\
&\qquad 
= \frac{1}{c} \sum_{i : t(i) \leq ct} \sum_{s = t(i)}^{t(i+1)-1}
    \mathbb{E}_{z \sim \mu^{\pi^\star_i}}\left[
                \sigma^\star_{r,i}\big(\bm Z(\infty)\big)
    \right]
    \\
    &\qquad = 
 \frac{1}{c} \sum_{s=cu}^{ct} \E_{\bm z\sim \mu^{\pi^\star}(s)} \sbrac{\sigma^\star_r (\bm z)}\, .
\end{align*}

So we have established that Term~\eqref{eqn:term1} tends to zero as $c \to \infty$ by applying the Azuma-Hoeffding inequality, and leveraging timescale separation, Term~\eqref {eqn:term2} also vanishes as $c \to \infty$. So term~\eqref{eqn:term3} is the only term remaining. 
Thus, we see that
\begin{align}
\bar{D}_r(t) 
&= \lim_{c\rightarrow \infty} \frac{1}{c} \sum_{i:t(i) \leq ct} \sum_{s=t(i)}^{t(i+1)-1} \E_{\bm z\sim \mu^{\pi^\star_i}} \sbrac{\sigma^\star_r (\bm z,\bar{\bm Q}^c(t(i)/c))}
\label{eq:Dlim2}    \\
&= \lim_{c\rightarrow \infty} \frac{1}{c} \sum_{s=cu}^{ct} \E_{\bm z\sim \mu^{\pi^\star}(s)} \sbrac{\sigma^\star_r (\bm z)} , .
\label{eq:Dlim1}
\end{align}
Thus, we see that \eqref{eq:Dlim} holds as required.  
\Endproof

Here, we see that due to timescale separation, the departure process on the fluid scale is ultimately determined by the stationary distribution of servers. This observation is critical both for characterizing the non-standard capacity region of this system and for finding throughput optimal policies.

We now see that the stationary distribution of WARP determines the limit of the departure process.  In the next section, we can use this to establish the optimality of the departure process.

\subsubsection{Optimality of Limit Points.}

So far, we have established that the limit of any convergent subsequence of $\brac{\bar{\bm Q}^c}_{c\in \Nats}$ satisfies fluid model equations~\eqref{eqn:fl_1}-\eqref{eqn:fl_2}, and we have seen that a timescale separation determines the departure process. We now need to verify \eqref{eqn:fl_3} with the following proposition. 

\begin{restatable}{proposition}{OptProp}\label{prop:opt} For any agnostic policy $\bar{\pi}\in \mathcal A$ the following inequality holds
  \begin{equation}
\label{eqn:mdp_optimal}
\int_s^t \sum_{r \in \mathcal R} \bar{ Q}_r(u) 
     d\bar{D}_r(u) \geq \int_s^t \sum_{r \in \mathcal R}\bar{ Q}_r(u)  \E_{\bm z \sim \mu^{\bar{\pi}}}[\sigma_r]du\, .
\end{equation}  
\end{restatable}

\Beginproof{}
We quickly note the following technical point before proceeding with the proof: $\bar{D}_r(u)$ is an increasing Lipschitz continuous function, and thus, its derivative exists everywhere except on a set of measure zero. Thus, we can take the Lebesgue integral of the derivative of $\bar{D}'_r(u)$ when it exists. See \cite{dudley2002real}, for instance, for further details. We also note that we continue to use the notation for $t(i)$ and $\pi^\star_i$ as defined in the proof of Proposition \ref{prop:timescale}.

    To see why \eqref{eqn:mdp_optimal} holds, we observe that the prelimit summations are Riemann integral approximations for the integrals above. In particular, notice that by \eqref{eq:Dlim}
\begin{align}\label{eq:d1}
    \int_s^t \bar{Q}_r(u) \, \mathbb{E}_{z \sim \mu^{\bar{\pi}}}[\sigma_r] \, du
    &= 
    \lim_{c \to \infty} \frac{1}{c} \sum_{u = \lfloor cs \rfloor + 1}^{\lfloor ct \rfloor}
        \bar{Q}_r^c\left( \left\lfloor u \right\rfloor_{\star} \right)
        \, \mathbb{E}_{z \sim \mu^{\bar{\pi}}}[\sigma_r]
    \\[2ex]
    \int_s^t \bar{Q}_r(u) \, d\bar{D}_r(u)
    &= 
    \lim_{c \to \infty} \frac{1}{c} \sum_{u = \lfloor cs \rfloor + 1}^{\lfloor ct \rfloor}
        \bar{Q}_r^c\left( \left\lfloor u \right\rfloor_{\star} \right)
        \, \mathbb{E}_{z \sim \mu^{\pi^\star_u}}
            [\sigma_r^\star]
    \label{eq:d2b}
\end{align}
where $\left\lfloor u \right\rfloor_{\star} := \max \{ t(i) : t(i) \leq t \} /c$
 is the most recent time index where the policy was updated, and we let $\pi^\star_u$ be the corresponding policy.

Notice that by the definition of $\pi^\star$ being the optimal average MDP solution, c.f. \eqref{eq:AREMDP}, we have for all $u$ 
\begin{equation}\label{eq:d3}
\E_{\bm z\sim \mu^{\pi^\star_u}} \Bigg[
\sum_{r \in \mathcal R}\bar{ Q}_r^c\big({{\left\lfloor u \right\rfloor_{\star} }}\big)  \sigma^\star_r
\Bigg]
\geq 
       \E_{\bm z \sim \mu^{\bar{\pi}}} \Bigg[
     \sum_{r \in \mathcal R}
     \bar{ Q}_r^c\big({\left\lfloor u \right\rfloor_{\star }} \big)\sigma_r
     \Bigg].
\end{equation}

Combining \eqref{eq:d1}, \eqref{eq:d2b} and \eqref{eq:d3} we see that  
\begin{align*}
        \int_s^t \sum_{r \in \mathcal R} \bar{ Q}_r(u) 
     \bar{ D}'_r(u) du
    &=
    \lim_{c\rightarrow\infty}\frac{1}{c} \sum_{u=\floor{cs}+1}^{\floor{ct}} \sum_{r \in \mathcal R}\bar{ Q}_r^c\big({{\left\lfloor u \right\rfloor_{\star} }}\big)  \E_{\bm z\sim \mu^{\pi^\star_u}} \sbrac{\sigma_r^\star(\bm z)}
    \\
    &=
    \lim_{c\rightarrow\infty} \frac{1}{c}\sum_{u=\floor{cs}+1}^{\floor{ct}}\sum_{r \in \mathcal R}
    \E_{\bm z\sim \mu^{\pi^\star_u}} \Big[
\sum_{r \in \mathcal R}\bar{ Q}_r^c\big({{\left\lfloor u \right\rfloor_{\star} }}
\big)  \sigma_r^\star(\bm z)
\Big]
\\
&
\geq 
\lim_{c\rightarrow\infty} \frac{1}{c}
\sum_{u=\floor{cs}+1}^{\floor{ct}}
\sum_{r \in \mathcal R}
       \E_{\bm z \sim \mu^{\bar{\pi}}}
       \Big[
     \sum_{r \in \mathcal R}
     \bar{ Q}_r^c\big({{\left\lfloor u \right\rfloor_{\star} }} \big)\sigma_r
     \Big]
     \\
     &
     =
         \lim_{c\rightarrow\infty} \frac{1}{c}\sum_{u =\floor{cs}+1}^{\floor{ct}} 
         \sum_{r \in \mathcal R}
         \bar{ Q}_r^c\big({{\left\lfloor u \right\rfloor_{\star} }} \big)\cdot  \E_{\bm z \sim \mu^{\bar{\pi}}}[\sigma_r(\bm z)]
    =
        \int_s^t \sum_{r \in \mathcal R}\bar{ Q}_r(u)  \E_{\bm z \sim \mu^{\bar{\pi}}}[\sigma_r(\bm z)]du.
\end{align*}
Thus, we see that~\eqref{eqn:fl_3} holds as required.
\Endproof

We can now conclude the proof of Theorem \ref{thm:fluid_convergence}. We have now shown, as required, the tightness of the sequence $\brac{\bm \lambda^c,\bar{\bm{Q}}^c,\bar{\bm{A}}^c,\bar{\bm{D}}^c, \bar{\bm Z}^c}_{c \in \mathbb{N}}$ and that any limit point of this sequence satisfies the fluid model equations \eqref{eqn:fl_1}, \eqref{eqn:fl_2} and \eqref{eqn:fl_3}. This completes the proof  of Theorem \ref{thm:fluid_convergence}.\hfill $\square$

\subsection{Proof of Theorem~\ref{thm:stochastic_stability}}
\label{append:stochastic_stability}

To prove the positive recurrence of the Markov process $\bar{\bm X}^c$ under the policy $\pi^\star$, Proposition~\ref{prop:Tightness}, Theorem~\ref{thm:fluid_stability} and general stability results from \cite{dai1995positive, bramson2008stability, dai1995stability} are applied.
The main idea is to leverage the $L_1$ convergence of a subsequence of $\brac{\bar{\bm Q}^c(t)}_{c \in \Nats}$ in conjunction with the multiplicative Foster's Lemma.

Note that if the policies $\pi^\star$ were not throughput optimal, then the following converse claim must hold:
\begin{equation}
\label{eq:A}
\parbox{\dimexpr\linewidth-4em}{
    \strut
        There exists an $\epsilon>0$ such that for all $c_\epsilon\in\mathbb N$ there exists a $c>c_\epsilon$ and an arrival rate vector $\bm \lambda^c \in \mathcal C^\epsilon$ such that the policy $\pi^\star$ is not positive recurrent. 
    \strut
  }
\end{equation}

We will now argue that this cannot be the case. 

We consider the sequence of arrival rates and queueing networks $\brac{\bm \lambda^c,\bar{\bm{Q}}^c,\bar{\bm{A}}^c,\bar{\bm{D}}^c, \bar{\bm Z}^c}_{c\in \Nats}$ as described in \eqref{eq:A} each operating under the policy $\pi^\star$.
Since $\mathcal C^\epsilon$ is compact, any limit point of $\bm \lambda^c$  belongs to $\mathcal C^\epsilon$.
Recall that a sequence of random variables that converges in distribution and is uniformly integrable, also converges in $L_1$.
By Proposition~\ref{prop:Tightness} and Theorem \ref{thm:fluid_convergence}, $\brac{\bm \lambda^c,\bar{\bm{Q}}^c,\bar{\bm{A}}^c,\bar{\bm{D}}^c, \bar{\bm Z}^c}_{c}$ has a subsequence that converges to $\brac{\bm \lambda,\bar{\bm{Q}},\bm \lambda t,\bar{\bm{D}}^c, \bm 0 t }$ where $\bar{\bm Q}$ fluid solution with arrival rate $\bm \lambda \in \mathcal C^\epsilon$. Let the subsequence be $\brac{\bm \lambda^{c_k},\bar{\bm{Q}}^{c_k},\bar{\bm{A}}^{c_k},\bar{\bm{D}}^{c_k}, \bar{\bm Z}^{c_k}}_{k \in \mathbb{N}}$. From the fluid stability Theorem~\ref{thm:fluid_stability} 
notice that $T$ can be chosen uniformly for all $\bm \lambda \in \mathcal C^\epsilon$. This observation will help us prove throughput optimality.
We know that there exists a $T>0$ such that for all $\bm \lambda \in \mathcal C^\epsilon$ $\bar{Q}_r(T)=0$ for $r\in \mathcal{R}$. Therefore, the subsequence $\brac{\bar{\bm Z}^{c_k}(T),\bar{\bm Q}^{c_k}(T)}_{k\in \Nats}$ converges in distribution to $(\bar{\bm Z}(T),\bar{\bm Q}(T))=(\bm 0,\bm 0)$.
Next, we claim that the subsequence $\brac{\bar{\bm Z}^{c_k}(T),\bar{\bm Q}^{c_k}(T)}_{k\in \Nats}$ is uniformly integrable.

It is important to note that the sequence $\brac{\bar{\bm Z}^c(t),\bar{\bm Q}^c(t)}_{c \in \Nats}$ is uniformly integrable for any $t\geq0$. This is because the unscaled server process ${\bm Z}^c(t)$ is uniformly bounded and because the queue length process $\bar{\bm Q}^c(t)$ is bounded above by the cumulative request arrival process $\bar{\bm A}^c(t)$ plus $\bar {\bm {Q}}^c(0)$. Moreover, we have $\mathbb{E}[(\bar{\bm A}^c(t))^2]$ is bounded as the increments of $\bar{\bm A}^c(t)$ are assumed to have bounded variance. Hence, ensuring uniform integrability of $\brac{\bar{\bm Q}^c(t)}_{c \in \Nats}$. 

Since the subsequence $\brac{\bar{\bm Z}^c(t),\bar{\bm Q}^c(t)}_{c \in \Nats}$ is uniformly integrable and $\brac{\bar{\bm Z}^c(T),\bar{\bm Q}^c(T)}$ converges in distribution to $(\bar{\bm Z}(T),\bar{\bm Q}(T)) = (\bm 0,\bm 0)$, we also have $L_1$ converges for the subsequence
\begin{equation}
\lim_{k\to \infty} \E \norm{\brac{\bar{\bm Z}^{c_k}(T),\bar{\bm Q}^{c_k}(T)}}_1=\E \norm{\brac{\bar{\bm Z}^{c_k}(T),\bar{\bm Q}^{c_k}(T)}}_1=0.
\end{equation}
The above $L_1$ convergence implies that there exists a $c_{\epsilon}$ such that for all $c>c_{\epsilon}$ we have
\begin{equation}
\label{eqn:ss_1}
\E \norm{\brac{\bar{\bm Z}^{c_k}(T),\bar{\bm Q}^{c_k}(T)}}_1 <(1-\delta),
\end{equation}
for any $\delta>0$.
Therefore, for $c=\norm{\brac{\bar{\bm Z}^{c_k}(T),\bar{\bm Q}^{c_k}(T)}}_1>c_\epsilon$ we can write 
\begin{equation}
\label{eqn:ss_2}
\E\sbrac{\norm{(\bm Z^c(cT) ,\bm Q^c\brac{cT})}_1-\norm{(\bm Z^c(0),\bm Q^c(0))}_1\mid \bm Q(0)}\leq -\delta c,
\end{equation}
where the inequality follows from~\eqref{eqn:ss_1}. The inequality in~\eqref{eqn:ss_2} shows that the conditions of the \textit{Multiplicative Foster's Lemma} are met [see Proposition 4.6 of~\cite{bramson2008stability}]. Therefore, the request queue process $\bar{\bm{Q}}^c$ under the policy $\pi^\star$ is positive recurrent for all $c > c_\epsilon$. Thus, we see that \eqref{eq:A} cannot hold. 
\hfill $\blacksquare$

\section{Azuma-Hoeffding Lemma}
We analyze a sequence of nested Martingale difference sequences using the following Azuma-Hoeffding Lemma. Arguments of this type are regularly used for mixing bounds in reinforcement learning \cite{qu2020finite}.


\begin{lemma}
\label{lem:azzuma_hoeffding}
Let $$m_r(i,s)= \sbrac{\sigma^\star_r\brac{\bm Z(t(i)+s),\bar{\bm Q}^c(t(i)/c)}-\E \sbrac{\sigma^\star_r\brac{\bm Z(t(i)+s),\bar{\bm Q}^c(t(i)/c)} \mid \mathcal{F}_{t(i)} }}, \ r\in \mathcal{R},$$ 
where $\mathcal{F}_{t(i)}$ is the filtration generated by $\left(\bm{X}(s) : 0 \leq s \leq t(i)\right)$\footnote{From the Markov property it implies that we can condition on filtration rather the process $\bm X$.} and $\left|m_r(i,s)\right|\leq B$ for all $i>0$, $s>0$.  If $\tau(c) \leq D c^{1-\epsilon}$ for some $\epsilon > 0$ and $D>0$, then for any $\delta>0$ we have the following bound
\begin{equation} \label{eq:lem1}
\mathbb{P}\left(
    \sup_{u \leq t} \left| 
        \frac{1}{c} 
        \sum_{i:t(i) \leq c u} 
        \sum_{s=0}^{t(i+1)-t(i)-1} m_r(i,s)
    \right| \geq \delta
\right)
\leq 2 D c^{1-\epsilon}\exp\left( -\frac{c^{2\epsilon} \delta^2}{2 B^2 t} \right).
\end{equation}
Therefore, we have that, with almost surely, uniformly on compact time intervals that
\[
\left| \frac{1}{c}
\sum_{i:t(i) \leq c u} 
\sum_{s=0}^{t(i+1)-t(i)-1}
    m_r(i,s)
\right|
\xrightarrow[c \rightarrow \infty]{} 0.
\]
\end{lemma}

\Beginproof{}
Using a union bound, we can write 
\begin{align*}
 \mathbb{P}\brac{\sup_{u\leq t}\left| \sum_{s=1}^{\tau(c)-1}\sum_{i:t(i) \leq cu} m_r(i,s)\right|
 \geq c\delta}
 &
 \leq\mathbb{P}\brac{\exists s\leq \tau(c)-1:\sup_{u\leq t} \left|\sum_{i:t(i) \leq c u} m_r(i,s)\right|\geq \frac{c\delta}{\tau(c)}}\\
  &
 \leq\mathbb{P}\brac{\exists s\leq D c^{1-\epsilon}:\sup_{u\leq t} \left|\sum_{i:t(i) \leq c u} m_r(i,s)\right|\geq \frac{c\delta}{D c^{1-\epsilon}}}\\
&
\leq \sum_{s=1}^{D c^{1-\epsilon}} \mathbb{P}\brac{\sup_{u\leq t}\left|\sum_{i:t(i) \leq cu} m_r(i,s)\right|\geq \frac{c\delta}{D c^{1-\epsilon}}}  \, . 
\end{align*}
Above, we note in the inequality that 
\[
t(i+1) - t(i) = \tau(Q(t(i))) \leq \tau(c+c \lambda_{\max}t) = D c^{1-\epsilon}.
\]
Then we note that $\sum_{i:t(i) \leq c} m_r(i,s)$ is a martingale difference sequence. Hence, the result follows by applying the Azuma-Hoeffding Inequality (See \cite{hoeffding1994probability}) to the last term in the above inequality. This can be applied to the maximum of the process by combining Azuma-Hoeffding with Doob's maximal inequality.

For almost sure convergence, we note that the summation of the \eqref{eq:lem1} if finite for all $\delta>0$. By the Borel-Cantelli Lemma, this implies almost convergence below $\delta$, and since this holds for all $\delta >0$, this implies almost sure convergence to zero as required.
 
\Endproof

\section{A Further Instability Counter-Example}
\label{append:instability_maxweight}

We have already given a counter-example demonstrating that MaxWeight is not throughput optimal with deterministic server lifetimes and arrival processes. We now extend this argument for Bernoulli arrivals and fixed decoherence probabilities.

\begin{restatable}{theorem}{MWInstable2}\label{thrm:MWInstable2}
    The MaxWeight policy is not optimal for throughput when servers have probabilistic decoherences.
\end{restatable}

\begin{figure}[h]
    \centering
    \includegraphics[width=0.5\textwidth, trim=0cm 2cm 0cm 
    5cm, clip]{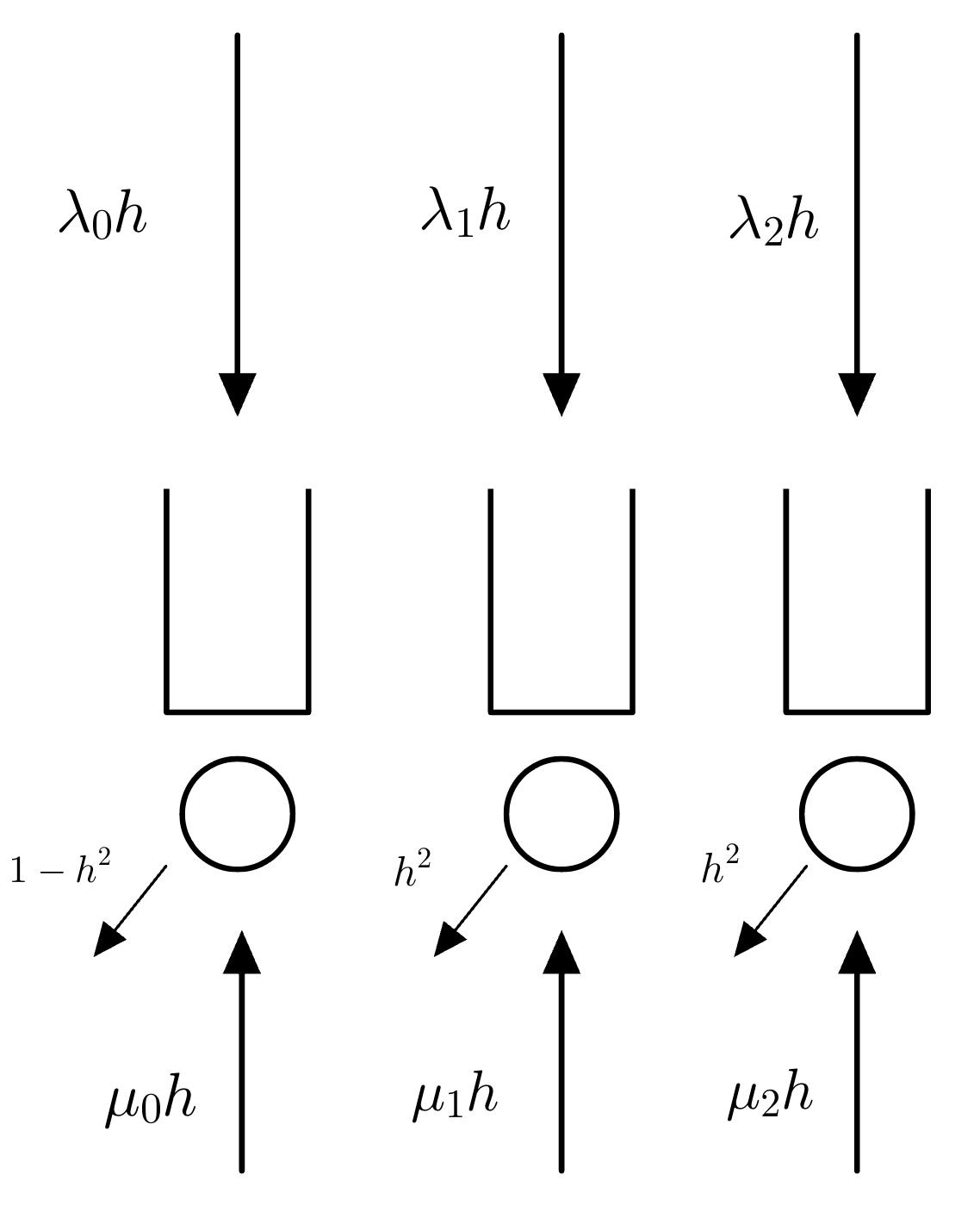}
    \caption{Queueing network topology with transitions}
    \label{fig:queue_network}
\end{figure}

\Beginproof{}
\noindent \textit{Outline.} We develop a counter-example for the network described in Figure \ref{fig:queue_network}. Here, all transitions are Bernoulli random variables, unlike our previous counterexample.
In particular, request arrivals occur as Bernoulli processes with probabilities $\lambda_{i}h $ for $i = 0, 1, 2$ servers arrive with probabilities $\mu_i h$ for $i = 0, 1, 2$. In both cases, $h$ is a small but positive number.
Queue $0$ has a high probability of decoherence, $1-h^2$, whereas Queues $1$ and $2$ have a low probability of decoherence, $h^2$. For the network depicted in Figure~\ref{fig:queue_network}, type $0$ requests are served jointly by type $0$, type $1$, and type $2$ servers; type $1$ requests are served by type $1$ servers; and type $2$ requests are served by type-$2$ servers.

We consider two policies: MaxWeight and the policy that prioritizes Queue $0$, which we call P$0$. Find settings where the network is unstable under MaxWeight but stable under the P$0$ policy. Thus demonstrating that MaxWeight cannot be throughput optimal.

\medskip
\noindent \textit{A Continuous Time Limit.} To briefly give some intuition for the proof, notice if we let $h \rightarrow 0$, and we rescale time into intervals of length $h$, then our network converges to a continuous time Markov chain with arrival rates $\lambda_i$ and service rates $\mu_i$ for $i = 0, 1, 2$, where servers at Queues $1$ and $2$ have infinite lifetimes and the server at Queue $0$ decoheres instantaneously, unless it is served at the moment it arrives. Given the unavailability of Queue 0 servers, it is not too hard to show that for this continuous-time Markov chain, MaxWeight will prioritize Queues $1$ and $2$, whereas by definition, P$0$ prioritizes Queue 0. Thus, for MaxWeight, we only serve Queue 0 when Queues $1$ and $2$ are empty, which occurs at a proportion of $(1-\lambda_1/\mu_1)(1-\lambda_2/\mu_2)$. (Note Queues 1 and 2 behave independently under MaxWeight and are individually empty with probability $(1-\lambda_i/\mu_i)$, $i=1,2$.) Similarly, for the P$0$ policy, we only serve Queues $1$ and $2$ at a proportion of $(1-\lambda_0/\tilde \mu_0)$. 
(Here $\tilde \mu^{-1}_0$ is the mean service rate of Queue 0, which is the time for two servers to arrive at Queue $1$ and $2$ and for a server to arrive at Queue $0$. The exact formula for this is given below in (B).)

This leads to the following condition for the instability of MaxWeight and two conditions for the stability of the priority policy: 
\begin{align}
    \mu_0\bigg(1 - \frac{\lambda_1}{\mu_1}\bigg)\bigg(1 - \frac{\lambda_2}{\mu_2}\bigg) &< \lambda_0, \tag{A} \label{eqn:instability_maxweight} \\
    \frac{1}{\tilde{\mu}_0} := \frac{1}{\mu_1} + \frac{1}{\mu_2} - \frac{1}{\mu_1+\mu_2} + \frac{1}{\mu_0} &< \frac{1}{\lambda_0}, \tag{B}\label{eqn:p0stability1} \\
    \bigg(1 - \frac{\lambda_0}{\tilde{\mu}_0}\bigg) > \frac{\lambda_1}{\mu_1}\, , \quad 
    \text{and} 
    \quad \bigg(1 - \frac{\lambda_0}{\tilde{\mu}_0}\bigg) &> \frac{\lambda_2}{\mu_2}.
    \tag{C}\label{eqn:p0stability2}
\end{align}
The first condition~\eqref{eqn:instability_maxweight} is the condition for the instability of Queue $0$ under MaxWeight. The second condition~\eqref{eqn:p0stability1} is the condition for the stability of Queue 0 under the P$0$ policy. (Note that service requires waiting for two servers to arrive at Queues $1$ and $2$ and for a server to arrive at Queue $0$. The maximum of two exponential random variables plus an additional exponential has a mean of $	\tilde \mu^{-1}_0$ defined as above. Thus, the condition ensures the interarrival times are longer than service times.) The third condition~\eqref{eqn:p0stability2} is the condition for the stability of Queues $1$ and $2$ under the P$0$ policy. (Note that the long-run service rate of Queue $1$, which occurs when Queue $0$ is idle, is $\mu_1 (1-\lambda_0/\tilde{\mu}_0)$, and to be stable, this rate must be greater than $\lambda_1$.)

Thus, if we can find parameters for which~\eqref{eqn:instability_maxweight}-\eqref{eqn:p0stability2} are satisfied, we see that P$0$ is stable, and MaxWeight is not. Thus, MaxWeight cannot be optimal for throughput. This covers the case of the continuous-time model.
Some related proofs for CTMCs can be found in the paper of~\cite{bonald2001impact}. 

Now it remains to identify the parameters for which the conditions in~\eqref{eqn:instability_maxweight}-\eqref{eqn:p0stability2} are satisfied and to verify that the DTMC satisfies these conditions for sufficiently small $h$.

\medskip
\noindent \textit{Parameters Satisfying~\eqref{eqn:instability_maxweight}-\eqref{eqn:p0stability2}.} We verify that there exist parameters for which~\eqref{eqn:instability_maxweight}-\eqref{eqn:p0stability2} are satisfied. In particular, we can take

\[
\lambda_1 = \lambda_2 = 150, \quad \mu_1=\mu_2 = 200, \quad \lambda_0 = 4, \quad \mu_0 = 20 \, .
\]
For Condition~\eqref{eqn:instability_maxweight}:
\[
    \mu_0\bigg(1 - \frac{\lambda_1}{\mu_1}\bigg)\bigg(1 - \frac{\lambda_2}{\mu_2}\bigg)
    = 20 \times \frac{1}{4} \times \frac{1}{4} = 1.25 < 4 = \lambda_0 \, .
\]
For Condition~\eqref{eqn:p0stability1}:
\[
\frac{1}{\tilde{\mu}_0} := \frac{1}{200} + \frac{1}{200} - \frac{1}{400} + \frac{1}{20} = \frac{23}{400}  < \frac{1}{4} = \frac{1}{\lambda_0} \, .
\]
For Condition~\eqref{eqn:p0stability2}:
\[
\bigg(1 - \frac{\lambda_0}{\tilde{\mu}_0}\bigg) = 1 - \frac{92}{400} = \frac{308}{400} > \frac{3}{4} = \frac{\lambda_1}{\mu_1}= \frac{\lambda_2}{\mu_2} \, .
\]
Thus, we see that there are parameters for which~\eqref{eqn:instability_maxweight}-\eqref{eqn:p0stability2} hold.

\medskip
\noindent \textit{Analysis for DTMC.} We have now completed the main conceptual components of the proof. What follows below is a detailed proof of the above inequalities, which accounts for the discrete transitions in the prelimit network, which depend on $h$. Here we show that~\eqref{eqn:instability_maxweight}-\eqref{eqn:p0stability2} hold for sufficiently small $h$. Moreover, we show that each expression is correct up to terms that are o(1) as $h \rightarrow 0$.\footnote{Such inequalities can likely be verified by a weak convergence argument. However, we provide direct proof here.}

Notice that excluding the transition for the decoherence of the type $0$ server, all transitions have a probability that is $O(h)$. Thus, the probability of a pair of any such transitions occurring simultaneously is $O(h^2)$. In particular, the probability of two transitions occurring simultaneously is bounded above by $h^2 C_0^{*}$ where $C_0^{*} = \max_{i=0,1,2}\{\mu^2_i, \lambda^2_i, 1\}$.  Thus, the probability of two transitions occurring simultaneously is $O(h^2)$ is bounded above by $ C_1^{*} = \binom{7}{2}  C_0^{*}$. Here, $\binom{7}{2}$ corresponds to the 28 pairs of transitions that can occur. We focus on verifying the instability condition for MaxWeight~\eqref{eqn:instability_maxweight} and then give conditions~\eqref{eqn:p0stability1} and~\eqref{eqn:p0stability2}.

\medskip
\noindent \textit{\eqref{eqn:instability_maxweight} Instability Condition for MaxWeight.}
We find parameters where Queue $0$ is starved out by the other two queues. We do this by lower bounding Queue $0$ and upper bounding Queues $1$ and $2$. In particular, we consider the lower-bound process where we assume every time there are two simultaneous transitions attempted, we remove a request from Queue $0$. 
This allows us to analyze situations where only one event occurs at a time since the worst that can happen from a pair of events is already factored into Queue $0$'s performance.

Now, notwithstanding the pairs of events just discussed, notice that whenever there is a server in Queue $1$ or $2$, that queue will be served if there is work to be done. In other words, unless two transitions occur simultaneously, Queues $1$ and $2$ have priority in effect over type $0$. Thus, excluding transitions with probability $O(h^2)$, Queue $0$ can only be served if both Queues $1$ and $2$ are empty.

\medskip
\noindent \textit{Idle Queue Conditions.}
Notice that for Queues $1$ (and Queue $2$), there is only one transition of order $O(h)$ that is influenced by other queues, 
that is, when the queue is empty, and there is service of a server waiting at Queue $1$. (This is because the server could be used to serve Queue $0$ or Queue $1$.)
This transition is always bounded below by $\lambda_1$, the probability of 
an arrival at Queue $1$. In other words, the transition at the queue and thus 
the queue size process at Queue $1$ can be bounded below by a single server 
queue with arrival probability $\lambda_1 h$ and departure probability $\mu_1 h - O(h^2)$. Thus, we can bound above the time that the Queue $1$ spends idle with a server waiting by:
\[
P(\text{Queue $1$ is idle with a server}) \leq 1 - \frac{\lambda_1}{\mu_1} - O(h)\, .
\]
We can apply the same coupling to Queue $2$, and since the processes considered are independent
we have that
\[
P(\text{Both Queue 1 \& 2 idle with a server}) \leq \bigg(1 - \frac{\lambda_1}{\mu_1}\bigg)\bigg(1 - \frac{\lambda_2}{\mu_2}\bigg) + O(h).
\]
Thus, we must wait for a server to arrive. Thus, the rate of departure can be bounded above by:
\[
h \mu_0 \bigg(1 - \frac{\lambda_1}{\mu_1}\bigg)\bigg(1 - \frac{\lambda_2}{\mu_2}\bigg) + O(h^2) \, .
\]
Thus, with the above upper-bound on the long run departure rate from Queue $0$, we arrive at the condition that Queue $0$ will be unstable for $h$ sufficiently small, provided $\lambda$ and $\mu$ satisfy:
\begin{equation}
       \mu_0 \bigg(1 - \frac{\lambda_1}{\mu_1}\bigg)\bigg(1 - \frac{\lambda_2}{\mu_2}\bigg) < \lambda_0 \, . 
\end{equation}

\medskip
\noindent \textit{\eqref{eqn:p0stability1} Mean Service Time.} 
Conversely, suppose we prioritize Queue $0$. That is, we only use servers for Queues $1$ and $2$ when Queue $0$ is empty. In that case, 
let us bound the service rate at Queue $0$. Notice that since the lifetime 
of server's at Queue $1$ and $2$ are $O\left({h^{-2}}\right)$, thus the probability 
of the only way of serving Queue $0$ (with probability that is not $O(h^2)$) is 
that both servers at Queue $1$ and $2$ arrive (and do not depart). Then a server 
arrives at Queue $0$ and is served.

Note the server arrival at Queues $1$ and $2$ is the maximum of two geometric random variables, 
with mean $\frac{1}{h\mu_1} + \frac{1}{h\mu_2} - \frac{1}{h\mu_1+h\mu_2} + O(1)$, and the time for the server arrival at Queue $0$ is 
geometric with mean $\frac{1}{h\mu_0}$. Thus, the mean time between service 
epochs at Queue $0$ is:
\[
\frac{1}{h\mu_1} + \frac{1}{h\mu_2} - \frac{1}{h\mu_1+h\mu_2} + \frac{1}{h\mu_0} + O(1) \,.
\]
In other words, we have that the following condition is sufficient for stability:
\[
\frac{1}{\tilde \mu_0} := \frac{1}{\mu_1} + \frac{1}{\mu_2} - \frac{1}{\mu_1+\mu_2} + \frac{1}{\mu_0} < \frac{1}{\lambda_0} . 
\]

That is, the mean time between service epochs at Queue $0$ is less than the mean time between arrivals at Queue $0$. This gives Condition~\eqref{eqn:p0stability1}.
(Notice that the service rate in (1) could be quicker if servers are waiting in advance at Queues $0$ and $1$. Thus,~\eqref{eqn:p0stability1} is a pessimistic and sufficient condition for Queue $0$ to be stable.)

\medskip
\noindent \textit{\eqref{eqn:p0stability2} Sufficient Condition for Stability.} From the above service rate, we see that the time that the Queue $0$ spends idle is at least:
\[
\left(1 - \frac{\lambda_0}{\tilde \mu_0}\right) \, .
\]
Once idle, Queue 1 requests are served at rate $h \mu_1$. We have the following sufficient 
condition for stability:
\[
\mu_1 \left(1 - \frac{\lambda_0}{\tilde \mu_0}\right) > \lambda_1 \, ,
\]
or equivalently:
\[
\left(1 - \frac{\lambda_0}{\tilde{\mu}_0}\right) > \frac{\lambda_1}{ \mu_1}\, . 
\]
This, along with the corresponding condition for Queue $2$, gives Condition~\eqref{eqn:p0stability2}.

\medskip
\noindent \textit{Conclusion.}
In summary, we see that conditions~\eqref{eqn:instability_maxweight}-\eqref{eqn:p0stability2} are satisfied for sufficiently small $h$. We have given parameters for which P$0$ can be stable but MaxWeight is not stable. Thus, MaxWeight is not optimal for throughput in the setting of an MDPN with entanglement memories.

\Endproof

\begin{rem}
We note that by adding further queues, Queues $3, 4, 5, ..., N$, we can include further products in condition~\eqref{eqn:instability_maxweight}. This will further exacerbate MaxWeight's instability. This suggests that the loss of throughput can decrease multiplicatively with the network size.
\end{rem}

   \end{document}